\documentclass[12pt]{article} 

\usepackage{defs1}
\usepackage{slashed}
  \usepackage{xcolor}
\title{Foliated manifolds, algebraic $K$-theory, and a secondary invariant}
\author{Ulrich Bunke\thanks{NWF I - Mathematik,
Universit{\"a}t Regensburg,
93040 Regensburg,
GERMANY, ulrich.bunke@mathematik.uni-regensburg.de}  
}

\theoremstyle{definition}

\newcommand{\alg}{\mathrm{alg}}
\newcommand{\triv}{\mathrm{triv}}
\newcommand{\per}{\mathrm{per}}
\newcommand{\op}{\mathrm{op}}
\newcommand{\filt}{\mathrm{filt}}
\renewcommand{\flat}{\mathrm{flat}}
\newcommand{\fr}{\mathrm{fr}}
\renewcommand{\max}{\mathrm{max}}
\newcommand{\Adams}{\mathrm{Adams}}

\newcommand{\an}{\mathrm{an}}

\renewcommand{\cl}{\mathrm{cl}}
\newcommand{\fol}{\mathrm{fol}}
\newcommand{\bbF}{\mathbf{F}}
\newcommand{\codim}{\mathrm{codim}}

\newcommand{\DGA}{\mathbf{CDGA}}
\newcommand{\KU}{\mathbf{KU}}

\renewcommand{\DD}{\mathbf{DD}}

\newcommand{\bM}{\mathbf{M}}

\newcommand{\bKU}{\mathbf{KU}}
\newcommand{\cK}{\mathcal{K}}

\newcommand{\CAlg}{{\mathbf{CAlg}}}

\renewcommand{\Proj}{\mathbf{Proj}}

\newcommand{\PSh}{{\mathbf{PSh}}}

\newcommand{\bK}{{\mathbf{K}}}

 \newcommand{\Vect}{{\mathbf{Vect}}}
 \newcommand{\CommGroup}{{\mathbf{CommGroup}}}
 \newcommand{\CommMon}{{\mathbf{CommMon}}}
 \newcommand{\Cat}{{\mathbf{Cat}}}

\newcommand{\reg}{{\tt reg}}

\newcommand{\ku}{{\mathbf{ku}}}

\renewcommand{\Dirac}{\slashed{D}}

\begin{document}
\maketitle
 \begin{abstract}
We introduce a $\C/\Z$-valued invariant of a foliated manifold  with a stable  framing  and with a partially flat vector bundle. This invariant can be expressed in terms of integration in differential $K$-theory, or alternatively, in terms of $\eta$-invariants of Dirac operators and local correction terms.  

Initially, the construction of the  element in $\C/\Z$ involves additional choices. But if  the codimension of the foliation is sufficiently small, then this element 
 is independent of these choices and therefore an invariant of the data listed above.

We  show that the invariant comprises various classical invariants like Adams' $e$-invariant, the $\rho$-invariant of twisted Dirac operators, or the Godbillon-Vey invariant from foliation theory.

Using methods from differential cohomology theory we construct a regulator map from the algebraic $K$-theory of smooth functions on a manifold   to its  connective $K$-theory with $\C/\Z$ coefficients. 
Our main result is a formula for the invariant in terms of this regulator and  integration in 
algebraic  and topological $K$-theory. 
 \end{abstract}
 \setcounter{tocdepth}{2}
\tableofcontents

\section{Introduction}

In this paper we introduce and analyse  an invariant 
$$\rho(M,\cF,\nabla^{I},s)\in \C/\Z$$
of an odd-dimensional closed spin manifold $M$ equipped with a real foliation $\cF$, a complex vector bundle with flat partial connection $\nabla^{I}$ in the direction of the foliation, and a stable framing $s$ of the foliation.
In order to define this number we must choose in addition a Riemannian metric on $M$, an extension of the partial connection $\nabla^{I}$ to a connection, and similarly, an extension of the canonical  flat partial  connection
on the normal bundle $\cF^{\perp}$ of the foliation.

Without any further conditions the number $\rho(M,\cF,\nabla^{I},s)$ may depend non-trivially on the additional geometric choices.  But
if the codimension of the foliation $\cF$  is sufficiently small, namely, if  \begin{equation}\label{fefwfwefewfewfewewf}
2\codim(\cF)< \dim(M)\ ,
\end{equation} then
$\rho(M,\cF,\nabla^{I},s)$ does not depend on the additional choices. 

\bigskip

The quickest way to define the invariant in Definition \ref{flwefjwefewff} is to use the integration in differential complex $K$-theory $\widehat{KU}^{*}$. Alternatively,  $\rho(M,\cF,\nabla^{I},s)$ can also be expressed as a combination of $\eta$-invariants of twisted Dirac operators and correction terms involving integrals of characteristic forms and their transgressions, see Proposition \ref{fjwfkljfklfjwelkfjkewjfewlkf9798237982749237432243}.

\bigskip

The invariant $\rho(M,\cF,\nabla^{I},s)$ is very interesting since it combines
various classical  secondary invariants in spectral geometry, topology and foliation theory in one object.
We will reveal these relations by analysing  special cases in Section \ref{jbjkhwekfewfwefewf89798}. We will observe that
 $\rho(M,\cF,\nabla^{I},s)$ subsumes Adams' $e$-invariant for framed manifolds, the rho-invariant for Dirac operators twisted with flat bundles,
and classical invariants from foliation theory like the Godbillon-Vey invariant.

\bigskip 

While the construction of the invariant $\rho(M,\cF,\nabla^{I},s)$ and the verification of its basic properties
are not very deep and based on  well-known methods from differential geometry and local index theory 
we think that its relation  with algebraic $K$-theory is much less obvious.  
 In the present paper we reveal this relation in the special case of 
 a foliated manifold of the form $$(M,\cF)=(P\times X,T_{\C}P\boxplus \{0\})\ .$$
  Here the manifold $P$ is closed and stably framed. A complex vector bundle
$(V,\nabla^{I})$ with flat partial connection on $M$ provides an algebraic $K$-theory class  of the ring of complex-valued smooth functions $C^{\infty}(X)$. We will write this class   as $$f^{o_{s}}_{!}([V,\nabla^{I}]^{\alg})\in K_{p}(C^{\infty}(X))$$  with notation to be introduced in Section  \ref{keklwfewfewfewf}, where
  $p:=\dim(P)$.  If $p >\dim(X)$ (this is exactly condition \eqref{fefwfwefewfewfewewf}), then we can  define a regulator transformation
$$\reg_{X}:K_{p}(C^{\infty}(X))\to \ku\C/\Z^{-p-1}(X)\ ,$$ see Definition \ref{klfwefewfewfwf}.
If we now assume that $X$ is a closed  spin manifold such that $p+\dim(X)$ is odd, then we have an integration map
$$\pi_{!}^{o}:\ku\C/\Z^{-p-1}(X)\to \ku\C/\Z^{-p-\dim(X)-1}(*)\cong \C/\Z\ ,$$ where
$\pi:X\to *$ and $o$ is the orientation of $\pi$ for $\ku\C/\Z$ induced by the spin structure.
Our main result is Theorem \ref{flkfefwefwefewfef}:  
\begin{theorem} \label{edhjlqkdjqwdqwdwqdw}
$$\rho(M,\cF,\nabla^{I},s)=\pi_{!}^{o}(\reg_{X}(f^{o_{s}}_{!}([V,\nabla^{I}]^{\alg})))\ .$$
\end{theorem}
The proof of this theorem will be finished in Section \ref{keklwfewfewfewf}. It is based on the diagram \eqref{dewewdedewdewd2342343} which comprises various Riemann-Roch type squares for integration in algebraic and topological $K$-theory and their differential refinements.


 \bigskip
 
In Section \ref{dkqwldqwdwqdwqdwqdwqd} we put the regulator $\reg_{X}$ into its natural general framework. We introduce the algebraic $K$-theory spectrum $\bK(M,\cF)$ of a foliated manifold and the Hodge-filtered connective complex $K$-theory spectrum $\ku^{\flat}(M,\cF)$. The regulator $\reg_{X}$ used in the theorem above is then a special case of a regulator
 $$\reg:\bK(M,\cF) \to \ku^{\flat}(M,\cF) \ .$$
 In order to justify to call this map  a regulator
  consider a complex manifold $M $ as a real manifold with a complex foliation $\cF:=T^{0,1}M$. In this case $\bK(M,T^{0,1}M)$  is the algebraic $K$-theory spectrum of $M$ defined using holomorphic vector bundles. Furthermore,  the homotopy groups of 
 $\ku^{\flat}(M,T^{0,1}M) $
  are the $\ku$-theory analogues of the integral Deligne cohomology groups.
 The regulator is an integral refinement of a version of Beilinson's regulator. 
 We will explain all this in detail in Section \ref{dkqwldqwdwqdwqdwqdwqd}.

 \bigskip
 
In Section \ref{gbiorgergergergreg}   we introduce   basic definitions from the theory of foliated manifolds and characteristic classes. 
The experienced reader could skip this section  in a first reading and use it as a reference for notation and normalization conventions. In Section \ref{jfkelwfwefewfewfewfewfwef} we give a quick introduction to the features of differential complex $K$-theory which are used in the construction of the invariant $\rho(M,\cF,\nabla^{I},s)$. The actual  construction of this invariant   will be given in  Section \ref{flkwjefwkejfeflkewjflkefjelkfjelfjewlfewf9790}. As mentioned above, in   Section  \ref{fkweflwefewfewfewfwfwfw} we provide
a spectral theoretic interpretation of the invariant, and in Section   \ref{jbjkhwekfewfwefewf89798}  relate it to various classical secondary invariants. 

In Section \ref{keklwfewfewfewf} we develop the theory which is necessary to state and prove Theorem \ref{edhjlqkdjqwdqwdwqdw}. Finally, Section  \ref{dkqwldqwdwqdwqdwqdwqd}  is devoted to the algebraic $K$-theory of foliated manifolds and the regulator in general. This section has a substantial overlap with the work of Karoubi \cite{karoubi45}, \cite{karoubiast} \cite{karoubi43}. In a certain sense it reformulates his constructions using the new technology  of the  $\infty$-categorical approach to $K$-theory and regulators developed in \cite{Bunke:2013ab}, \cite{Bunke:2012fk}, \cite{Bunke:2013aa}, \cite{Bunke:2014aa}.

\bigskip

{\em Acknowledgement:  This work was partially supported by the SFB 1085 ''Higher Invariants'' of the DFG.}

\section{Basic notions}\label{gbiorgergergergreg}

\subsection{Foliated manifolds}\label{fewl453534535435}

We introduce the category of foliated manifolds $\Mf_{\C\textrm{-}\fol}$ and its full subcategory $\Mf_{\fol}$ of manifolds with real foliations. 

\bigskip

In the present paper we consider foliations from the infinitesimal point of view. 
Textbook references for real foliations (with  emphasis on the decomposition of the manifold into leaves) are
\cite{goodbillon-book} or \cite{schochetmoore}. We further refer to  \cite{tondeur} which  also contains a comprehensive list of further  references for foliation theory in general.  For the
  complex foliation associated to a complex structure on a manifold 
   we refer to the textbook sections \cite[IX.2]{MR1393941} and \cite[Ch. 1.3]{MR2359489}.

\bigskip

Let $M$ be a smooth manifold and let $T_{\C}M:=TM\otimes_{\R}\C$ be the complexified tangent bundle. A section  of $   T_{\C}M$ is called a complex vector field. A complex vector field $X\in \Gamma(M,T_{\C}M)$ acts as a derivation $(X,f)\mapsto X(f)$ on the algebra $C^{\infty}(M)$ of complex-valued smooth functions $f$.   For a pair of complex vector fields  $X,Y $ we can consider the  commutator $[X,Y]\in \Gamma(M,T_{\C}M)$. It is the unique complex vector field such that $$[X,Y](f)=X(Y(f))-Y(X(f))$$
for all $f\in C^{\infty}(M)$.

\bigskip

If $\cF\subseteq T_{\C}M$ is a subbundle, then we have an inclusion $\Gamma(M,\cF)\subseteq \Gamma(M,T_{\C}M)$ of spaces of sections.
 \begin{ddd} A subbundle $\cF \subseteq T_{\C}M$ is called integrable if for any two sections $X,Y\in \Gamma(M,\cF )$ we also have $[X,Y]\in \Gamma(M,\cF)$.
\end{ddd}

\begin{ddd}
A  foliation of a smooth manifold $M$ is an integrable subbundle $\cF  \subseteq T_{\C}M$. A foliated manifold is a  pair $(M,\cF )$ of a manifold and a   foliation.  \end{ddd} 

Since $T_{\C}M$ is the complexification of  the real vector bundle $TM$ we have a complex antilinear  involution $X\mapsto \bar X$. For a subbundle $\cF\subseteq T_{\C}M$ we let $\bar \cF\subseteq T_{\C}M$ denote the subbundle obtained by applying this automorphism to the elements of $\cF$.
\begin{ddd} 
A foliation is called real if $ \bar{\cF}=\cF$.
   In this case we define the real  integrable  subbundle $\cF_{\R}:=\cF\cap TM\subseteq TM$. 
   \end{ddd}

\bigskip

Let $f:M\to N$ be a smooth map between manifolds. Its  differential is a map of bundles $$df :T_{\C}M\to f^{*}T_{\C}N$$ over $M$.

\begin{ddd}
We say that $f:(M,\cF )\to (M^{\prime},\cF^{\prime} )$ is a foliated map if its differential preserves the foliations in the sense that
$df (\cF )\subseteq f^{*}\cF^{\prime} $.
\end{ddd}

The composition of two foliated maps is again a foliated map.

\begin{ddd}
We let $\Mf_{\C\textrm{-}\fol}$ denote the category of foliated manifolds and foliated maps.
We further let $\Mf_{\fol}\subset \Mf_{\C\textrm{-} \fol}$ be the full subcategory of foliated manifolds with real foliations.
\end{ddd}

Let $\Mf$ denote the category of smooth manifolds. Then we have functors
$$\Mf_{\fol}\to \Mf_{\C\textrm{-}\fol}\to \Mf\ ,$$
where the first is the  inclusion of a full subcategory, and the second forgets the foliation.
The category of foliated manifolds has a cartesian product. It is given by 
$$(M,\cF)\times (M^{\prime},\cF^{\prime})\cong (M\times M^{\prime},\cF\boxplus \cF^{\prime})\ .$$

\begin{ex}\label{groijoergerregr}{\rm 
If $M$ is a complex manifold, then the subbundle  $\cF :=T^{0,1}M\subseteq T_{\C}M$ is a 
complex foliation. Vice versa, a complex foliation $\cF $ with the additional property that $\cF \oplus \bar{\cF }\cong T_\C M$  equips  $M$ with a complex structure such that $T^{0,1}M=\cF$.
Moreover, a foliated map between such foliated manifolds is the same as a holomorphic map. 
}
\end{ex}

\begin{ex}\label{lkjwqdwdqwdqwdqwd}{\rm 
Every manifold $M$ has a minimal foliation $\cF_{\min}:=\{0\}$ and a maximal foliation $\cF_{\max} :=T_{\C}M$.
These foliations are real. If $M$ is equipped with the minimal foliation and $(M^{\prime},\cF^{\prime})$ is a foliated manifold, then every smooth map $M\to M^{\prime}$ is foliated. 
Similarly, if  $M^{\prime}$ is equipped with the maximal foliation, then for a foliated manifold $(M,\cF)$ every smooth map $M\to M^{\prime}$   is foliated.
}
\end{ex}

 \begin{ex}\label{djqlwdqwdqwdq}{\rm 
 Let $\pi:W\to B$  be a submersion. Then the complexification of the vertical bundle $T^{v}\pi:=\ker(d\pi)\subseteq TW$ 
 defines a real foliation $\cF^{v}$ called the vertical foliation. The map $\pi$ is foliated for any choice of a foliation on $B$.
 }
 \end{ex}

 \begin{ex}{\rm 
 Let $\Gamma$ be a discrete group which acts freely and properly on a manifold $\tilde B$ from the right with  quotient $B:=\tilde B/\Gamma$. Furthermore, let $X$ be a manifold with a left action of $\Gamma$. Then we consider the manifold
 $M:=\tilde B\times_{\Gamma} X$.  The vertical foliations $\tilde \cF^{v} $ and $\tilde \cF^{H} $ associated to the projections $\tilde B\times X\to \tilde B$ and $\tilde B\times X\to X$     descend to the quotient and define  the     vertical and horizontal foliations    $\cF^{v} $ and $\cF^{H} $  on  $M$. Note that $ \cF^{v} $ is the vertical foliation of the submersion $ M\to B$. We have $\cF^{H}\oplus \cF^{v}\cong T_{\C} M$.  }
 \end{ex}

 \begin{ex}\label{dekldlqwdqwdqwd}{\rm 

 Let $(B,\cF)$ be a foliated manifold. We call a map $f:W\to B$ transversal to $\cF$ if for every $w\in W$ we have the relation $\cF_{f(w)}+df(TW_{x})=T_{f(x)}B$.  If $f$ is transversal to $\cF$, then we can define
 a maximal foliation $f^{-1}\cF$ on $W$ such that $f$ becomes a foliated map. We must set $f^{-1}\cF:=df^{-1}(f^{*}\cF)$. 
 
 \bigskip

In particular, if $P$ is a manifold, then we can consider the projection $\pi:P\times B\to B$. 
 In this case, $ \pi^{-1}\cF=T_{\C}P\boxplus \cF$.   
 
 }
 \end{ex}

\begin{rem}{\rm 
In the main part of the present paper from Section \ref{ergregregerger345345} to Section \ref{keklwfewfewfewf} we will be concerned with real foliations. The more general case of complex foliations will be relevant in Section \ref{dkqwldqwdwqdwqdwqdwqd} where we introduce the algebraic $K$-theory of foliations and the regulator. The general notion of a complex foliation interpolates between the more classical real foliations on the one end, and the complex foliation associated to a complex structure on a complex manifold on the other end. 
}
\end{rem}

\subsection{Filtrations on the de Rham complex}\label{uieerggrgerg34}

A foliation on a manifold induces a decreasing multiplicative filtration of the de Rham complex, see e.g. \cite[1.1]{karoubi43}, \cite[Sec. 4]{tondeur}.

\bigskip

We consider a foliated manifold $(M,\cF)$. 
By $(\Omega(M),d)$ we denote the complexified de Rham complex of $M$.

\begin{ddd}\label{ilfjewlfwfewfewfewfwfw}
For $n,p\in \nat$ we define the subspace $$F^{p}\Omega^{n}(M)\subseteq \Omega^{n}(M)$$   of forms which vanish after the insertion of $n-p+1$ sections of $\cF$.
\end{ddd} 
Roughly speaking, a form belongs to $F^{p}\Omega^{n}(M)$ if it has at least $p$ legs in the direction normal  to the foliation.
The family of these subspaces for all $p$ forms a decreasing filtration of $\Omega^{n}(M)$. More precisely we have the following chain of inclusions
$$ \Omega^{n}(M)=F^{0}\Omega^{n}(M)  \supseteq F^{1}\Omega^{n}(M)\supseteq  \dots \supseteq F^{\codim(\cF)}\Omega^{n}(M)\supseteq F^{\codim(\cF)+1}\Omega^{n}(M)=0\ .$$ 
Combining these filtrations for all $n$ together we get 
a decreasing filtration $(F^{p}\Omega(M))_{p\in \nat}$ of the graded commutative algebra $\Omega(M)$ which  is multiplicative, i.e., the wedge product restricts to maps \begin{equation}\label{hjkqhdkdhwqkdwqd98789}
\wedge :F^{p}\Omega^{m}(M)\otimes F^{q}\Omega^{n}(M)\to F^{p+q}\Omega^{m+n}(M)\ .\end{equation}

These properties are in fact true for arbitrary subbundles $\cF$ of $T_{\C}M$. But as  a consequence of the integrability of $\cF$ this filtration is also  preserved by the de Rham differential, i.e.,
$(F^{p}\Omega(M),d)$ is a subcomplex of $(\Omega(M),d)$  for every $p\in \nat$.

\bigskip

\begin{ddd}
For   a foliated manifold $(M,\cF)$ we write $\Omega(M,\cF)$  for the de Rham complex $\Omega(M)$  considered as a filtered commutative differential graded algebra.
\end{ddd}

\begin{ex}{\rm If $(M,\cF)$ is a complex manifold (Example \ref{groijoergerregr}), then the filtration on $\Omega(M,T^{0,1}M)$
is called the Hodge filtration. 
}\end{ex}


If $f:(M,\cF)\to (M^{\prime},\cF^{\prime})$ is a foliated map, then $f^{*}:\Omega(M^{\prime},\cF^{\prime})\to \Omega(M,\cF)$ is a morphism of filtered commutative differential graded algebras.

\bigskip

We let $\DGA$ and $\DGA^{\filt}$ denote the categories of graded commutative differential graded algebras and filtered graded commutative differential graded algebras. For  categories $\bC,\bD$ we can consider the functor category  $  \Fun(\bC^{\op},\bD) $\footnote{Later  in the paper we will also use the notation 
$\PSh_{\bD}(\bC):= \Fun(\bC^{\op},\bD)$ and call these objects $\bD$-valued presheaves on $\bC$}.
 We can formalize the properties of the filtered de Rham complex discussed above by saying that we have  a functor  
 $$\Omega\in \Fun(\Mf_{\C\textrm{-}\fol}^{\op},   \DGA^{\filt}) \ .$$

\begin{ddd}\label{jkdjlqwdqwdqwd}   We define the  functors    $$DD^{-}\in \Fun(\Mf_{\C\textrm{-}\fol}^{\op},\DGA)\ , \quad DD^{\per}\in \Fun(\Mf^{\op}, \DGA)$$ 
  by $$DD^{-}(M,\cF)  :=\prod_{p\in \Z} F^{p}\Omega (M)[2p]\ , \quad DD^{\per}:=\prod_{p\in \Z} \Omega(M)[2p]\ .$$
\end{ddd}
We call $DD^{\per}$ the periodic and $DD^{-}$ the negative de Rham complex.
Note that $DD^{-}$ has a decomposition into a product of components $$DD^{-}(M,\cF)\cong \prod_{p\in \Z} DD^{-}(M,\cF)(p)\ , \quad   DD^{-}(M,\cF)(p) :=F^{p}\Omega (M)[2p]\ ,$$
and the product on $DD^{-}$ is induced   by the wedge products of forms \eqref{hjkqhdkdhwqkdwqd98789} componentwise as $$DD^{-}(M,\cF)(p)\otimes DD^{-}(M,\cF)(q)\to DD^{-}(M,\cF)(p+q)\ .$$
 The description of the product for $DD^{\per}$ is similar.

\begin{rem}{\rm
The cohomology of $DD^{\per}(M)$ is the two-periodic de Rham cohomology of $M$.
It is the natural target of the Chern character from topological $K$-theory, see Definition \ref{ffwefwefewfewfwfw}. 
The complex $DD^{-}(M,\cF)$ will receive characteristic forms for vector bundles with connections which are flat in the direction of the foliation, see Definition \ref{qldjqwldqwdqwdqwd}.
}\end{rem}

\subsection{Vector bundles with flat partial connections}\label{fhfjlwefkjfewfewfewfwf}

We introduce the notion of a vector bundle with a flat partial connection  on a foliated manifold.
While a connection (see e.g. \cite[Ch. III]{MR1393940}) allows to differentiate sections in all directions of the manifold, a partial connection only differentiates in the direction of the foliation. The typical example of a flat partial connection is a holomorphic structure $\bar \partial$ on a complex vector bundle on a complex manifold \cite[III.2]{MR2359489}.

\bigskip

We consider a foliated manifold $(M,\cF)$. 
Let $V\to M$ be a complex vector bundle. 

\begin{ddd}\label{lwfwfewfewfewf}
A partial connection on $V$ is a map
$$\nabla:\Gamma(M,V)\to \Gamma(M,\cF^{*}\otimes V)$$ which satisfies the  Leibniz rule.
 \end{ddd}

\begin{rem}{\rm
For   sections $X\in \Gamma(M,\cF)$ and $\phi\in \Gamma(M,V)$ we write as usual $$\nabla_{X}\phi:=i_{X}(\nabla(\phi))
\in \Gamma(M,V)$$ for the evaluation of $\nabla\phi$ at $X$. With this notation the Leibniz rule has the form
$$\nabla_{X}(f\phi)=X(f)\phi+f\nabla_{X}\phi\ , \quad \forall f\in C^{\infty}(M)\ , \quad   \forall\phi\in \Gamma(M,V)\ , \quad  \forall X\in \Gamma(M,\cF)\ .$$ 
}\end{rem}
 
The foliation gives rise to a graded commutative differential graded algebra whose underlying commutative graded algebra is given by  $\Omega(\cF):=\Gamma(M,\Lambda^{*}\cF)$. Its differential $d^{\cF}$  is fixed by the prescription  \begin{equation}\label{jdgjdhjdhk1j2hdkjhd2d2ud2du1d12d}
d^{\cF}:\Omega^{0}(\cF)\to \Omega^{1}(\cF)\ , \quad  d^{\cF}(\phi):=d\phi_{|\cF}\ ,
\end{equation}
where $d$ is the usual de Rham differential and we use the identification  $\Omega^{0}(\cF)= C^{\infty}(M)=\Omega^{0}(M)$.
 We further  write $\Omega(\cF,V):=\Gamma(M,\Lambda^{*}\cF\otimes V)$.
 As in the case of usual connections we can extend $\nabla$ uniquely  to a derivation on the $\Omega(\cF)$-module $\Omega(\cF,V)$.
Its curvature defined by $$R^{\nabla}:=\nabla^{2}\in \End( \Omega(\cF,V))$$ is $\Omega(\cF) $-linear and hence a two-form on $\cF$ with values in $\End(V)$, i.e., we have 
  $R^{\nabla}\in  \Omega^{2}(\cF, \End(V))$.
\begin{ddd}\label{ijflkewfewfoieoiuoiuoiuoiuwef}
A partial connection $\nabla$  on $V$ is called flat if $R^{\nabla}=0$.
\end{ddd}

We now consider a foliated map $f:(M^{\prime},\cF^{\prime})\to (M,\cF)$.
If $V\to M$ is a vector bundle with a partial connection $\nabla$, then $f^{*}V$ 
has an induced partial connection $f^{*}\nabla$.
It is characterized by
\begin{equation}\label{gf433g34g3353454353453tr}
\nabla_{X}(f^{*}\phi)=f^{*}\nabla_{df(X)}\phi\ , \quad \forall m^{\prime}\in M^{\prime}\ , \quad  \forall X\in T_{m^{\prime}}M^{\prime}\ , \quad  \phi\in \Gamma(M,V)\ .
\end{equation}
This formula has to be understood as an equality between elements in the fibre $(f^{*}V)_{m^{\prime}}$.
Because of the relation
$$f^{*}R^{\nabla}=R^{f^{*}\nabla}$$ 
the pull-back of a flat partial connection is again flat.

\begin{ex}{\rm If $M$ is a complex manifold with foliation $\cF=T^{0,1}M$, then a  flat partial connection on a complex vector bundle $V$ is the same as a holomorphic structure. In this situation the flat partial  connection   is usually denoted by $\bar \partial$.
}
\end{ex}

\begin{ex}{\rm 
If $M$ is equipped with the minimal foliation, then a  partial connection on a vector bundle is no additional data.
In the opposite case, where $M$ has the maximal foliation,   a flat partial connection is the same as a flat connection.}
\end{ex}
 \begin{ex}\label{fklwefjwefewf}{\rm 
Let $$\cF^{\perp}:=T_{\C}M/\cF$$ be the normal bundle of a foliation.
Then $\cF^{\perp}$ has a natural flat partial connection $\nabla^{\cF^{\perp}}$. It is given by
$$\nabla_{X}[Y]:= [\ [X,Y]\ ]\ ,$$ where $  X\in \Gamma(M,\cF)$ and the vector field $Y\in \Gamma(M,T_{\C}M)$ represents the section $[\ Y\ ] \in \Gamma(M,\cF^{\perp})$ of the normal bundle.}
\end{ex}

\begin{ex}\label{jfwlefjlwefewfewfewfw}{\rm 
Let $\pi:W\to B$ be a submersion and consider the vertical foliation $T^{v}\pi$ on $W$, see Example \ref{djqlwdqwdqwdq}.
If $V\to B$ is any vector bundle, then $\pi^{*}V\to W$ has a canonical flat partial connection $\nabla^{I}=\pi^{*}\nabla$, where $\nabla$ is the canonical flat partial connection
on $V$ in the direction of the trivial foliation. In view of \eqref{gf433g34g3353454353453tr}
it is characterized by the condition that for  $\phi\in \Gamma(B,V)$ we have  $\nabla^{I}\pi^{*}\phi=0$.
}
\end{ex}



%

\subsection{Connections and characteristic forms}

We introduce the Chern character forms and Chern forms of complex vector bundles with connection (see e.g.  the textbook sections \cite[Sec. 1.5]{MR2273508} and \cite[Sec. III.3]{MR2359489}). In the    foliated case we discuss the consequences of the fact that the connection extends a flat partial connection. These consequences are expressed in terms of the position of the characteristic forms with respect to the filtration  of the de Rham complex introduced in Section \ref{uieerggrgerg34}, and they are the classical starting point for the construction of characteristic  classes for foliations, see e.g.  \cite{MR0307250}, \cite{MR0402773},  \cite{MR0334237}, \cite{karoubi43}.

\bigskip

Let $(M,\cF)$ be a foliated manifold and let $(V,\nabla^{I})$ be a complex vector bundle with a flat partial connection.

\begin{ddd} \label{fjwelfwfewf42343242wf}
A connection $\nabla$ on  $V$ is  an extension of $\nabla^{I}$, if the relation
$\nabla_{X}\phi=\nabla^{I}_{X}\phi$ holds  for all $\phi\in \Gamma(M,V)$ and $X\in \Gamma(M,\cF)$.
\end{ddd}

One can show that a flat partial  connection admits extensions. Furthermore, the set of extensions of a flat partial connection is a torsor over the complex vector space
$$\Gamma(M, \cF^{\perp,*}\otimes \End(V))\ .$$

\begin{ex}{\rm 
A connection on $\cF^{\perp}$ which extends the  flat partial  connection $\nabla^{I,\cF^{\perp}}$ of Example \ref{fklwefjwefewf} is called a Bott connection.}
\end{ex}

 \begin{ex}\label{fjewflewjflewfewf}{\rm 
Let $f:W\to B$ be a submersion and let $\cF^{v}$ be the vertical foliation (Example \ref{djqlwdqwdqwdq}). If $V\to B$ is a complex vector bundle, then
$f^{*}V\to W$ has a canonical flat partial connection $\nabla^{I}$, see Example \ref{jfwlefjlwefewfewfewfw}. If $\nabla$ is any connection on $V$, then $f^{*}\nabla$ extends $\nabla^{I}$.

More generally, if $f:W\to B$ is transverse to a foliation $\cF$ on $B$ and $(V,\nabla^{I})$ is a vector bundle with flat partial connection on $(B,\cF)$, then $(f^{*}V,f^{*}\nabla^{I})$ is a vector bundle with flat partial connection on $(W,f^{-1}\cF)$, see Example \ref{dekldlqwdqwdqwd}. If $\nabla$ is a connection on $V$ extending $\nabla^{I}$, then $f^{*}\nabla $ is a connection on $f^{*}V$ extending $f^{*}\nabla^{I}$.

}
\end{ex}

If $\nabla$ is a connection on a complex vector bundle, then  we consider its curvature
$$R^{\nabla}:=\nabla^{2}\in \Omega^{2}(M,\End(V))\ .$$ The Chern character  form of $\nabla$ is the closed  inhomogeneous complex-valued form
$$\ch_{0}(\nabla)+\ch_{2}(\nabla)+\ch_{4}(\nabla)+\dots:=\Tr \exp\left(-\frac{R^{\nabla}}{2\pi i}\right)$$
with homogeneous components $\ch_{2p}(\nabla)\in \Omega^{2p}_{\cl}(M)$. We will consider the Chern character form as a zero cycle in the periodic complex $DD^{\per}(M)$.

\begin{ddd}\label{ffwefwefewfewfwfw}
We define
$$\ch(\nabla):=(\ch_{2p}(\nabla))_{p\in \Z}\in  Z^{0}(DD^{\per}(M))\ .$$
\end{ddd}

\begin{rem}\label{dkjqwlkdjqwlkjdlwqkdwdwqdqwd}{\rm  
In this remark we explain how the Chern character form  behaves under complex conjugation and inserting adjoint connections. First of all, 
if we choose a hermitean metric $h$ on the bundle $V$ with connection $\nabla$, then we can form the adjoint connection $\nabla^{*}$, which is characterized by the following relation in $\Omega^{1}(M)$:
$$dh(\phi,\psi)=h(\nabla\phi,\psi)+h(\phi,\nabla^{*}\psi)\ , \quad  \mbox{for all}\:\phi,\psi\in \Gamma(M,V)\ .$$ Applying $d$ to this equality again we get 
$$0=h(R^{\nabla}\phi,\psi)+h(\phi, R^{\nabla^{*}}\psi)\ .$$ In view of the $2\pi i$-factor in the definition of the Chern character form this equality implies  the relation \begin{equation}\label{ccewcwecweewwerewr}
\overline{\ch(\nabla)}=\ch(\nabla^{*})\ .
\end{equation} 
The connection $\nabla$ is called unitary if $\nabla^{*}=\nabla$. In this case, the Chern character form $\ch(\nabla)$ is real.}
\end{rem}

Let $(V,\nabla^{I})$ be a complex vector bundle on $(M,\cF)$ with a flat partial connection $\nabla^{I}$. We further choose a hermitean metric $h$ on $V$.

\begin{ddd}\label{kdlqwdqwdqwdqd}
A flat partial connection $\nabla^{I}$ is called unitary (with respect to $h$), if  $$ d^{\cF}h(\phi,\psi)=h(\nabla^{I}\phi,\psi)+h(\phi,\nabla^{I}\psi)$$ for all $\phi,\psi\in \Gamma(M,V)$.
\end{ddd}
See \eqref{jdgjdhjdhk1j2hdkjhd2d2ud2du1d12d} for $d^{\cF}$.

\begin{lem}\label{dhqwkdqkwddqwdwqdioipopioipoopi}
If $\nabla^{I}$ is unitary (with respect to $h$), then it admits a  unitary  extension $\nabla$. 
\end{lem}
\proof
 Let $\nabla_{0}$ be some extension of $\nabla^{I}$. Then $\nabla^{*}_{0}$ is a second extension of $\nabla^{I}$ and
$$\nabla:=\frac{1}{2}(\nabla_{0}+\nabla_{0}^{*})$$
is a unitary extension of $\nabla^{I}$. \hB

\bigskip

There is a  filtration on $\Omega(M,\End(V))$  defined similarly as the filtration of $\Omega(M,\cF)$ introduced in Definition \ref{ilfjewlfwfewfewfewfwfw}. It  is compatible with the $\Omega(M,\cF)$-module structure. 

\begin{lem}If $\nabla$ extends a flat partial connection $\nabla^{I}$ on $V$, then
$R^{\nabla}\in F^{1}\Omega^{2}(M,\End(V))$.
\end{lem}
\proof
We have
$$R^{\nabla}_{|\Lambda^{2}\cF }=R^{\nabla^{I}}=0\ .$$ \hB
 This lemma has consequences for the Chern character forms.
 
 \begin{kor}\label{fwefwefewfewfewfewf2434ewfewfwefewf}
 If $\nabla$ extends a flat partial connection on $V$, then
 $$\ch_{2p}(\nabla)\in F^{p}\Omega^{2p}_{\cl}(M,\cF)\ .$$
 \end{kor}

 \begin{ddd}\label{qldjqwldqwdqwdqwd}
  If $\nabla$ extends a flat partial connection on $V$, 
then we define  
$$\ch^{-}(\nabla):=(\ch_{2p}(\nabla))_{p\in \Z}\in  Z^{0}(DD^{-}(M,\cF))\ .$$
\end{ddd}

We let $\Vect^{\flat,\nabla}(M,\cF)$ and $\Vect^{\nabla}(M)$ denote the symmetric monoidal categories (with respect to the direct sum) of pairs $(V,\nabla)$ of complex vector bundles with connection, where in the first case $\nabla$ extends a flat partial connection. In both cases morphisms are connection preserving vector bundle morphisms.
 
 If $f:M^{\prime} \to M $ is a smooth map and $(V,\nabla)\in \Vect^{\nabla}(M)$, then 
we can define $(f^{*}V,f^{*}\nabla)\in \Vect^{\nabla}(M^{\prime})$ and have the relation \begin{equation}\label{wqdqwdqwdwqdwqwqdqd}
 f^{*}\ch(\nabla)=\ch(f^{*}\nabla)\ .\end{equation}
Similarly, if $f:(M^{\prime},\cF^{\prime})\to (M,\cF)$ is a foliated map and $(V,\nabla)\in \Vect^{\flat,\nabla}(M,\cF)$, then
$(f^{*}V,f^{*}\nabla)\in  \Vect^{\flat,\nabla}(M^{\prime},\cF^{\prime})$ and we have the relation
\begin{equation}\label{wqdqwdqwdwqdwqwqdqd1}
\ch^{-}(f^{*}\nabla)=f^{*}\ch^{-}(\nabla)\ .\end{equation}
Thus the Chern character forms are characteristic forms. In addition, they are additive, i.e., the Chern character form of a direct sum is the sum of the Chern character  forms of the summands.  
These properties will be important for the construction of the regulator in Subsection \ref{kfjwelfewfewfewfewfewfewfe}.

 

\bigskip

Let $(V,\nabla)$ be a complex vector bundle with connection on a manifold $M$. Then we define the Chern forms $$c_{p}(\nabla)\in \Omega^{2p}(M)$$ of $\nabla$ as the homogeneous components of the following inhomogenous form
$$1-c_{1}(\nabla)+c_{2}(\nabla)-\dots=\det\left(1-\frac{R^{\nabla}}{2\pi i}\right)\ .$$
The Chern forms can be expressed as homogeneous polynomials in the Chern character forms. In particular, if $(M,\cF)$ is a foliated manifold and $\nabla$ extends a flat partial connection, then we have \begin{equation}\label{jlkjlwejflkewjfewf89798}
c_{p}(\nabla)\in F^{p}\Omega^{2p}_{\cl}(M,\cF)\ .\end{equation}

\subsection{Characteristic forms of real foliations}

We introduce characteristic forms of real vector bundles on real foliated manifolds. We in particular discuss the $\hA$-form (see e.g.  \cite[Sec. 1.5]{MR2273508}).

\bigskip 

 We consider a real foliated manifold $(M,\cF)$. If $V\to M$ is a real vector bundle on a real foliated manifold $(M,\cF)$, then in analogy with Definitions \ref{lwfwfewfewfewf} and \ref{ijflkewfewfoieoiuoiuoiuoiuwef} we have the notion of a flat partial connection $$\nabla^{I}:\Gamma(M,V)\to \Gamma(M,\cF_{\R}^{*}\otimes V)$$ on $V$.
We furthermore have the notion of a connection $\nabla$ on $V$ extending $\nabla^{I}$ (compare with
Definition  \ref{fjwelfwfewf42343242wf}).

%

 \bigskip

Let $V\to M$ be a real vector bundle with connection $\nabla$. We let  $\nabla_{\C}$ denote  the induced connection on the complexification $V\otimes \C$.  
We   define the Pontrjagin forms of $\nabla$ by $$p_{k}(\nabla):=(-1)^{k}c_{2k}(\nabla_{\C})\in \Omega_{\cl}^{4k}(M)\ .$$ 
If $ \cF $ is a real  foliation and  $\nabla$ extends a flat partial connection, then by \eqref{jlkjlwejflkewjfewf89798} we have \begin{equation}\label{wwejoijiofwef9789}
p_{k}(\nabla)\in  F^{2k}\Omega^{4k}_{\cl}(M,\cF)\ .\end{equation}

\bigskip

In order to define the $\hA$-form we consider the  symmetric polynomials $p_{i}$ of degree $4i$ in variables $x_{\ell}$ of degree $2$ defined by the relation 
$$1+p_{1} +p_{2} +\dots=\prod_{\ell} (1-x_{\ell}^{2})\ .$$
We   define homogeneous polynomials $\hA_{4k} $ for $k\ge 1$ in the variables $p_{i}$  by the relation
$$ 1+\hA_{4} +\hA_{8} +\dots=\prod_{\ell}\frac{\frac{x_{\ell}}{2}}{\sinh(\frac{x_{\ell}}{2})}\ .$$
Then the components of the $\hA$-form are defined by
$$\hA_{4k}(\nabla):=\hA_{4k}(p_{1}(\nabla),p_{2}(\nabla),\dots)  \in \Omega_{\cl}^{4k}(M)\ .$$ 
 If $(M,\cF)$ is real foliated and the connection $\nabla$ extends a flat partial connection, then by \eqref{wwejoijiofwef9789} we have
$$\hA_{4k}(\nabla)\in F^{2k}\Omega^{4k}_{\cl}(M,\cF)$$

\bigskip

We again want to consider the $\hA$-form as a zero cycle of $DD^{\per}(M)$, or of $DD^{-}(M,\cF)$ in the foliated case. In order to simplify the notation we set     $\hA_{2k}(\nabla):=0$ if $k$ is odd.

\begin{ddd} If $\nabla$ is a connection on a real vector bundle on $M$, then we define
$$\hA(\nabla):=(\hA_{2p}(\nabla))_{p\in \Z}\in Z^{0}(DD^{\per}(M))\ .$$
If $(M,\cF)$ is a real foliated manifold and $\nabla$ extends a flat partial connection, then we define
$$\hA^{-}(\nabla):=(\hA_{2p}(\nabla))_{p\in \Z}\in Z^{0}(DD^{-}(M,\cF))\ .$$
\end{ddd}

The $\hA$-form is multiplicative, i.e., the $\hA$-form of a direct  sum of connections is the product of the $\hA$-forms of the summands. Furthermore, the $\hA$-form of a trivial connection is the multiplicative unit.

%
 
\begin{ex}\label{kjlwefwefwef}{\rm 
If $(M,\cF)$ is a real foliated manifold, then the real normal bundle $$\cF_{\R}^{\perp}:=TM/\cF_{\R}$$ of the foliation has a flat partial connection $\nabla^{I,\cF_{\R}^{\perp}}$ similar as in Example \ref{fklwefjwefewf}. If we choose a connection$ \nabla^{\cF_{\R}^{\perp}}$ extending $\nabla^{I,\cF_{\R}^{\perp}}$, then we obtain a cycle
$$\hA^{-}(\nabla^{\cF_{\R}^{\perp}})\in  Z^{0}(DD^{-}(M,\cF))\ .$$
}
\end{ex}

 \subsection{Transgression}

 We introduce the transgression (see \cite[Sec. 1.5]{MR2273508}) of characteristic forms and  discuss its basic properties. In the case of foliations we discuss
 the consequences of the fact that the connections extend a fixed flat partial connection (see e.g. \cite[Sec. 4]{karoubi43}). In this section we also introduce the notion of a stable framing and the conventions for the associated trivial connection.

 \bigskip
 
 We consider the unit interval $I:=[0,1]$ with coordinate $t$.
 For $i=0,1$ let  $\iota_{i}:*\to  I$ be the inclusions of the endpoints of the interval.
 Let $M$ be a smooth manifold and let $V\to M$ be a vector bundle. 
 Given two connections $\nabla_{0}$ and $\nabla_{1}$ we can consider a connection $\tilde \nabla$ on $\pi^{*}V\to I\times M$ such that $(\iota_{i}\times \id_{M})^{*}\tilde \nabla=\nabla_{i}$ for $i=0,1$.  For example we could take the linear interpolation $t\pi^{*}\nabla_{1}+(1-t)\pi^{*}\nabla_{0}$.  
 \bigskip
 
 The integration of forms along the fibre of $\pi:I\times M\to M$ is a map of graded vector spaces
 $$\int_{I\times M/M}:\Omega(I\times M)\to \Omega(M)[-1]\ .$$ It induces a map 
 $$\int_{I\times M/M} :DD^{\per}(I\times M)\to DD^{\per}(M)[-1]\ .$$ Since the interval $I$ has a non-empty boundary the integration is not a morphism of complexes. In fact,
 by Stokes' theorem we have the relation
 \begin{equation}\label{dqwdwqdwqdqdw}
(\iota_{1}\times \id_{M})^{*}-(\iota_{0}\times \id_{M})^{*}=d\circ \int_{I\times M/M}+\int_{I\times M/M}\circ d\ .
\end{equation} \begin{ddd}
 The transgression of the Chern character form is defined by
 $$\widetilde{\ch}(\nabla_{1},\nabla_{0}):=\int_{I\times M/M} \ch(\tilde \nabla)\in DD^{\per}(M)^{-1}/\im(d)\ .$$
 \end{ddd}
    From \eqref{dqwdwqdwqdqdw} and the facts that the Chern character forms are closed  and natural (see \eqref{wqdqwdqwdwqdwqwqdqd})
     we immediately get  the identity \begin{equation}\label{hfjwkjehfkjwehfewfewf897987}
 d \widetilde{\ch}(\nabla_{1},\nabla_{0}):=\ch(\nabla_{1})-\ch(\nabla_{0})\ .
\end{equation}
 One can check that  the transgression is independent  of the choice of the connection $\tilde \nabla$ interpolating between $\nabla_{0}$
 and $\nabla_{1}$. At this point it is relevant that we consider the transgression as a class modulo exact forms.
 Furthermore, we have the identities \begin{equation}\label{fewfwefwefwefew32434234234}
\widetilde{\ch}(\nabla_{1},\nabla_{0})+\widetilde{\ch}(\nabla_{2},\nabla_{1})+\widetilde{\ch}(\nabla_{0},\nabla_{2})=0\ , \quad \widetilde{\ch}(\nabla_{1},\nabla_{0})+\widetilde{\ch}(\nabla_{0},\nabla_{1})=0\ .  
\end{equation}
In order to see e.g. the first equality in \eqref{fewfwefwefwefew32434234234} one can integrate the Chern form of an interpolation  between
the three connections  over a two-simplex. 
 \begin{rem}{\rm 
 If we choose a hermitean metric on $V$, then we can form the adjoint connections. From   
  \eqref{ccewcwecweewwerewr} we get the relation \begin{equation}\label{qwdqwdqwdwqdqwdwqdqwdqwd}
\overline{\widetilde{\ch}(\nabla_{1},\nabla_{0})} =\widetilde{\ch}(\nabla^{*}_{1},\nabla^{*}_{0})\ .
\end{equation}
  }\end{rem}

\bigskip 
 
We now assume that $(M,\cF)$ is foliated and that the connections $\nabla_{i}$  for $i=0,1$ extend the same flat partial connection $\nabla^{I}$.
Then we can equip $I\times M$ with the foliation  $T_{\C}I\boxplus \cF$ introduced in Example \ref{dekldlqwdqwdqwd}. 
We can furthermore find an interpolation $\tilde \nabla$ which extends the flat partial connection $\pi^{*}\nabla^{I}$, e.g. the linear interpolation.

\bigskip

We now observe that the integration preserves the filtration, i.e., that it induces a map
$$\int_{I\times M/M}:F^{p}\Omega^{k}(I\times M,T_{\C}I\boxplus  \cF)\to F^{p}\Omega^{k-1}(M,\cF) \ .$$
At this point it is crucial that we included the tangent bundle of the interval into the foliation.
Hence we get an induced map
$$\int_{I\times M/M}:DD^{-}(I\times M,T_{\C}I\boxplus  \cF)\to DD^{-}(M,\cF)[-1]\ .$$
 
  \begin{ddd} Let  $(M,\cF)$ be a foliated manifold and let $(V,\nabla^{I})$ be a complex vector bundle with a flat partial connection on $M$. If $\nabla_{0}$ and $\nabla_{1}$ are two connections on $V$ extending $\nabla^{I}$, then we define 
 the transgression of the Chern character form   by
 $$\widetilde{\ch^{-}}(\nabla_{1},\nabla_{0}):=\int_{I\times M/M} \ch^{-}(\tilde \nabla)\in DD^{-}(M,\cF)^{-1}/\im(d)\ .$$
 \end{ddd}
 
Note again that $\widetilde{\ch^{-}}(\nabla_{1},\nabla_{0})$  is independent of the choice of the interpolation $\tilde \nabla$.
 
\begin{ex}{\rm We consider a foliated manifold $(M,\cF)$ and a complex vector bundle $V \to M$.
  If $p\in \nat$ is such that  $p>\codim(\cF)$,   then we have $F^{p}\Omega^{2p}(M)=0$.
  
  \bigskip

Assume that $\nabla^{I}_{0}$ and $\nabla^{I}_{1}$ are two flat partial connections on a complex vector bundle $V$ and let $\nabla_{0}$, $\nabla_{1}$ be corresponding extensions. If $p>\codim(\cF)$, then
  $\widetilde{\ch}_{2p}(\nabla_{1},\nabla_{0})$ is closed since by \eqref{dqwdwqdwqdqdw} and Corollary \ref{fwefwefewfewfewfewf2434ewfewfwefewf} its differential belongs to $F^{p}\Omega^{2p}(M)=0$.
  Its cohomology class does not depend on the choice of the extensions $\nabla_{1}$ and $\nabla_{0}$.
We therefore get a secondary characteristic class $$c(\nabla_{1}^{I},\nabla_{0}^{I}):=\widetilde{\ch}_{2p}(\nabla_{1},\nabla_{0})\in H^{2p-1}(M;\C)\ .$$
}\end{ex}

\begin{ex}{\rm
The following construction generalizes the Kamber-Tondeur  classes (introduced in this form in \cite{MR1303026})  to the foliated case.
Let $\nabla^{I}$ be a flat partial connection. If we choose a hermitean metric $h^{V}$ on $V$, then similarly as in the case of connections (see  Remark \ref{dkjqwlkdjqwlkjdlwqkdwdwqdqwd}) we can define
an adjoint flat partial connection $\nabla^{I,*}$.
It is  characterized by the relation
$$d^{\cF}h(\phi,\psi)=h(\nabla^{I}\phi,\psi)+h(\phi,\nabla^{I,*}\psi)\ , \quad  \phi,\psi\in \Gamma_{c}(M,V)\ .$$ 
 
 Let $\nabla$ be an extension of $\nabla^{I}$. Then $\nabla^{*}$ extends $\nabla^{I,*}$.
 We consider the  form $$\widetilde{\ch}_{2p}(\nabla,\nabla^{*})\in F^{p}\Omega^{2p-1}(M,\cF)/\im(d)\ .$$
By \eqref{qwdqwdqwdwqdqwdwqdqwdqwd} and \eqref{fewfwefwefwefew32434234234} we have the relation
 \begin{equation}\label{dewdwedwedewded}
\overline{\widetilde{\ch}_{2p}(\nabla,\nabla^{*})}= \widetilde{\ch}_{2p}(\nabla^{*},\nabla)= -\widetilde{\ch}_{2p}(\nabla,\nabla^{*})\ , 
\end{equation} i.e., the  form  $\widetilde{\ch}_{2p}(\nabla,\nabla^{*})$ is imaginary.
 
 \bigskip
 
  For $p>\codim(\cF)$ the class
$$c_{2p-1}(\nabla^{I}):=c(\nabla^{I},\nabla^{I,*})\in H^{2p-1}(M;\C) $$
does not depend on the choice of the hermitean  metric. 
By \eqref{dewdwedwedewded} it is imaginary, i.e., it belongs to the real subspace   $iH^{2p-1}(M;\R)\subseteq  H^{2p-1}(M;\C) $.
\bigskip

We can apply this construction to the bundle $\cF^{\perp}$ with its canonical flat partial  connection $\nabla^{I,\cF^{\perp}}$, see Example \ref{fklwefjwefewf}. The class \begin{equation}\label{kfkwejwlkefjlwekfjewfewfopipoi234}
c_{2p-1}(\nabla^{I,\cF^{\perp}})\in H^{2p-1}(M;\C)
\end{equation}   is closely related to the Godbillon-Vey class of the foliation.  If the foliation $\cF$ is real, then we can explain the place of this  invariant in the classification of characteristic classes for foliations. See Remark \ref{ergegojerglerogergeg}, in particular \eqref{ewfwefewfwf432342344123}.

}
\end{ex}

Let $V\to M$ be a real vector bundle. Then using a similar notation as above we can define
\begin{equation}\label{nkdqlwdwqddwqdj09809}
\widetilde{\hA}(\nabla_{1},\nabla_{0}):=\int_{I\times M/M}\hA(\tilde \nabla)\in DD^{\per}(M)^{-1}/\im(d)\ .\end{equation}
We have

 \begin{equation}\label{1e1h2ekj12ej12ke21hek2eh}d \widetilde{\hA}(\nabla_{1},\nabla_{0}):=\hA(\nabla_{1})-\hA(\nabla_{0})\ .\end{equation} The transgression is independent  of the choice of the connection $\tilde \nabla$.   Furthermore, we have the identities  $$\widetilde{\hA}(\nabla_{1},\nabla_{0})+\widetilde{\hA}(\nabla_{2},\nabla_{1})+\widetilde{\hA}(\nabla_{0},\nabla_{2})=0\ , \quad \widetilde{\hA}(\nabla_{1},\nabla_{0})+\widetilde{\hA}(\nabla_{0},\nabla_{1})=0\ . $$  
 
 Finally, if $\nabla_{1}$ and $\nabla_{0}$ extend the same flat partial connection, then we can define
 $$\widetilde{\hA^{-}}(\nabla_{1},\nabla_{0}):=\int_{I\times M/M}\hA(\tilde \nabla)\in DD^{-}(M)^{-1}/\im(d)\ .$$

Let $V\to M$ be a real vector bundle. By $V\oplus\underline{ \R^{n}}$ we denote the sum of $V$ with the trivial vector bundle of dimension $n$.
\begin{ddd}\label{hffhweuifhiuwefwfewfwefew}
A stable framing of $V$ is a trivialization of $V\oplus \underline{ \R^{n}}$ for some $n\in \nat$.  
\end{ddd}
A trivialization of $V\oplus \underline{\R^{n}}$ induces a  trivialization of $V\oplus \underline{\R^{n+k}}$ in a natural way. We will consider the corresponding two stable framings 
as equivalent. In particular, we can always arrange that two stable framings of $V$ trivialize the same bundle. 
A stable framing as above naturally induces a flat connection on $V\oplus\underline{ \R^{n}}$ which will be called the associated stable trivial connection and denoted by  $\nabla^{V,\triv}$ or $\nabla^{s}$. If $\nabla^{V}$ is another connection on $V$ and we want to consider the transgression of a characteristic  form between $\nabla^{V}$ and $\nabla^{V,\triv}$, then we secretly extend
$\nabla^{V}$ to $V\oplus \underline{ \R^{n}}$ using the trivial connection on the second summand.
The transgression will be invariant under further stabilization.

\bigskip

We now come back to the foliated manifold $(M,\cF)$ and   assume that the foliation $\cF$ is real.
Let $g^{TM}$ be a Riemannian metric on $M$.  
We get a decomposition $TM\cong \cF_{\R}\oplus \cF_{\R}^{\perp}$ and  choose   a connection $\nabla^{\cF_{\R}^{\perp}}$ extending $\nabla^{I,\cF_{\R}^{\perp}}$.  We  further assume that  $\cF_{\R}$ has a stable framing $s $.  Let $\nabla^{\cF_{\R},\triv}$ be the associated stable trivial connection on $\cF_{\R}$. We have  forms  
$$\widetilde{\hA}(\nabla^{LC},\nabla^{\cF_{\R},\triv}\oplus \nabla^{\cF_{\R}^{\perp}})\in DD^{\per}(M)^{-1}/\im(d) $$ and $$ \hA(\nabla^{\cF_{\R},\triv}\oplus \nabla^{\cF_{\R}^{\perp}})=\hA(\nabla^{\cF_{\R}^{\perp}})\in Z^{0}(DD^{-}(M))\ .$$

\section{Differential $K$-theory}\label{jfkelwfwefewfewfewfewfwef}

In this subsection we recall some basic features of the Hopkins-Singer version of differential complex $K$-theory.
 
 \bigskip
 
%
References for the following material are the foundational paper by Hopkins-Singer \cite{MR2192936}, but also \cite{MR2732065} and \cite{MR2664467}. For differential orientations and Umkehrmaps we refer to \cite{MR2664467}, \cite{Freed-Lott}, and \cite{2012arXiv1208.3961B}

\subsection{Basic structures}

We describe the differential extension $\widehat{KU}^{*}$ of the generalized cohomology theory $\KU^{*}$. Here $    \KU^{*}$ is the periodic topological complex $K$-theory which is represented by the spectrum $\KU$.
For every $p\in \Z$ we have a 
 contravariant functor
$$\widehat{KU}^{p}:\Mf^{\op}\longrightarrow  \Ab $$
from smooth manifolds to abelian groups.  This functor is connected with periodic topological complex $K$-theory $\KU^{*}$ via a transformation $$I:\widehat{KU}^{p}\to \KU^{p}$$ of abelian group valued functors. The transformation $I$ 
 maps differential $K$-theory classes to their underlying topological $K$-theory classes.
 Furthermore, differential $K$-theory is connected with differential forms through natural transformations $R$ and $a$.
The curvature $R$ is a natural transformation
$$ R:\widehat{KU}^{p}\to Z^{p}(DD^{\per}(M))  \ .$$
  In particular, if $x\in \widehat{KU}^{p}(M)$, then $R(x)\in Z^{p}(DD^{\per}(M))$
is a differential form representing the Chern character of the underlying topological class $I(x)$.
The  transformation $$a:  DD^{\per,p-1}/\im(d) \to \widehat{KU}^{p}$$ encodes the secondary information contained in differential $K$-theory classes. 
     All these structures and their compatibilities are nicely encoded in the following commutative diagram, also called the  differential cohomology diagram \cite{MR2365651}:
  $$\xymatrix{& DD^{\per,p-1} /\im(d)\ar[dr]^{a}\ar[rr]^{d}&&Z^{p}(DD^{\per}  )\ar[dr]^{\mathrm{Rham}}&\\H^{p-1}(DD^{\per} ) \ar[ur]\ar[dr]&&\widehat{KU}^{p} \ar[ur]_{R}\ar[dr]^{I}&&H^{p}(DD^{\per} )\ . \\&\KU\C/\Z^{p-1} \ar[ur]\ar[rr]^{\mathrm{Bockstein}}&&\KU^{p} \ar[ur]^{\ch}&}$$
 Its upper and lower parts are segments of long exact sequences, and the diagonals are exact at the center.
 
 \bigskip 
 
The flat part of differential $K$-theory  is defined as the kernel of the curvature transformation $R$: $$\widehat{KU}^{p}_{\flat}:=\ker(R:\widehat{KU}^{p}\to Z^{p}(DD^{\per}))\ .$$ It is canonically isomorphic to  K-theory with coefficients in $\C/\Z$ (with a shift):\begin{equation}\label{eq700hhh}\widehat{KU}^{p}_{\flat} \cong \KU \C/\Z^{p-1} \ .\end{equation}
Here $ \KU \C/\Z$ is an abbreviation for the spectrum $\KU\wedge \mathbf{M}\C/\Z$, where $\mathbf{M}\C/\Z$ is the Moore spectrum for the abelian group $\C/\Z$.

 \bigskip

The sequence
\begin{equation}\label{eq700}\KU^{p-1}\stackrel{\ch}{\to} DD^{\per,p-1}/\im(d)\stackrel{a}{\to} \widehat{KU}^{p}\stackrel{I}{\to} \KU^{p}\to 0\end{equation}
is exact.

 \bigskip

The differential $K$ -theory of a point is given by 
\begin{equation}\label{djhqkjwdqwdq9879}\widehat{KU}^{p}(*)\cong \left\{\begin{array}{cc} \Z&\mbox{if $p$ is even}\\
\C/\Z&\mbox{if $p$ is odd\ .}\end{array}\right.\end{equation}

 \bigskip

Differential $K$-theory is not homotopy invariant. The deviation from homotopy invariance is quantified by the homotopy formula. If $\hat x\in \widehat{KU}^{p}([0,1]\times M)$, then the homotopy formula states that
\begin{equation}\label{homotopyformula}(\iota_{1}\times \id_{M})^{*}\hat x-(\iota_{0}\times \id_{M})^{*}\hat x=a\left(\int_{[0,1]\times M/M} R(\hat x)\right)\ . \end{equation}

\subsection{The cycle map}

A complex vector bundle with connection $(V,\nabla)$  on a manifold $M$ gives rise to a differential $K$-theory class
$[V,\nabla]\in \widehat{KU}^{0}(M)$ such that  \begin{equation}\label{dkdqwdwqdqwdwqdkl8997}R([V,\nabla])=\ch(\nabla)\in Z^{0}(DD^{\per}(M))\ , \quad I([V,\nabla])=[V]\in \KU^{0}(M) \ .\end{equation} 

Let $M\mapsto \pi_{0}(\Vect^{\nabla}(M))$ denote the contravariant  functor which associates to the manifold $M$ the commutative monoid (induced by the direct sum) of isomorphism classes of pairs $(V,\nabla)$ of vector bundles with connection on $M$ and to a smooth map between manifolds $f:M\to M^{\prime}$ the pull-back $f^{*}$.
 The additive natural transformation
$$[\dots]:\pi_{0}(\Vect^{\nabla})\to \widehat{KU}^{0}$$ is called the cycle map and fits into the commuting diagram  of natural transformations between monoid-valued functors
$$\xymatrix{&Z^{0}(DD^{\per})\\\pi_{0}(\Vect^{\nabla})\ar[r]^{[\dots]}\ar[ur]^{[V,\nabla]\mapsto \ch(\nabla)}\ar[dr]_{[V,\nabla]\mapsto[V]}&\widehat{KU}^{0}\ar[u]^{R}\ar[d]_{I}\\&\KU^{0}}\ .$$

For compact manifolds $M$ the cycle map is known to be surjective \cite{MR2732065}. Assume that $\nabla_{0}$ and $ \nabla_{1}$ are two connections on the same complex vector  bundle $V$. Then as a consequence of the homotopy formula   \eqref{homotopyformula} we have the equality 
\begin{equation}\label{kdqklwdjqwlkdjqlwdkwqd89798}[V,\nabla_{1}]-[V,\nabla_{0}]=a(\tilde \ch(\nabla_{1},\nabla_{0}))\end{equation}
in $\widehat{KU}^{0}(M)$.
\subsection{Differential orientations and Umkehr maps}\label{wlekfjwelfjewlfwef123}

Let $\pi:W\to B$ be a proper submersion such that the vertical bundle $T^{v}\pi$ has a $\mathrm{Spin}^{c}$-structure. Then $\pi$ is  equipped with an orientation  $o$ (called the Atiyah-Bott-Shapiro orientation \cite{MR0167985}) for the cohomology theory $\KU^{*}$ and admits an Umkehr map
$$\pi^{o}_{!}:\KU^{p}(W)\to \KU^{p-d}(B)\ ,$$ where $d:=\dim(W)-\dim(B)$ is the dimension of the fibre (assume for simplicity that $B$ is connected). Since $\KU\C/\Z$ is a $\KU$-module spectrum we also have an integration
\begin{equation}\label{hjkfhjkhkjhkjefwefefwe7987}
\pi^{o}_{!}:\KU\C/\Z^{p}(W)\to \KU\C/\Z^{p-d}(B)\ .
\end{equation}

 The $\KU$-orientation $o$ determines a cohomology class
$$\hA(o)\in H^{0}(DD^{\per}(W))$$ such that the Riemann-Roch theorem holds true:
$$\xymatrix{\KU^{p}(W)\ar[d]^{\pi_{!}^{o}}\ar[r]^{\ch}&H^{p}(DD^{\per}(W))\ar[d]^{\int_{W/B}\hA(o)\cup \dots}\\ \KU^{p-d}(B)\ar[r]^{\ch}&H^{p-d}(DD^{\per}(B))}$$

Differential refinements of  $\KU$-orientations have been studied in detail in \cite{MR2664467}, \cite{Freed-Lott}, see also 
\cite{2012arXiv1208.3961B} for a more homotopy-theoretic approach.
In order to refine the $\KU$-orientation $o$ to a  $\widehat{KU}$-orientation $\hat o$ we must choose additional geometric structures. First of all we choose a metric on the vertical tangent bundle $T^{v}\pi$ and a horizontal distribution $T^{h}\pi$, i.e., a complement of the vertical bundle in $TW$. These structures induce a vertical Levi-Civita connection $\nabla^{T^{v}\pi}$, see e.g \cite[Prop. 10.2]{MR2273508}. In order to fix $\hat o$ we must further choose a $\mathrm{Spin}^{c}$-extension $\tilde \nabla^{T^{v}\pi}$ of $\nabla^{T^{v}\pi}$. The $\widehat{KU}$-orientation $\hat o$
 gives rise to an Umkehr map (see \cite[3.2.3]{MR2664467})
$$\pi^{\hat o}_{!}:\widehat{KU}^{p}(W)\to \widehat{KU}^{p-d}(B)\ .$$ The $\widehat{KU}$-orientation further provides a form $$\hA(\hat o)\in Z^{0}(DD^{\per}(W))$$ representing the class $\hA(o)$.  The Umkehr map in  $\widehat{KU}$-theory fits into the  following commutative diagram. \begin{equation}\label{r23r23r23r23r23r235435346546}
\hspace{-1cm}\xymatrix{\KU\C/\Z^{p-1}(W)\ar[d]^{\pi_{!}^{o}}\ar@/^1cm/[rr]&DD^{\per}(W)^{p-1}/\im(d)\ar[d]^{\int_{W/B}\hA(\hat o)\wedge \dots}\ar[r]^(0.6){a}&\widehat{KU}^{p}(W)\ar[d]^{\pi_{!}^{\hat o}}\ar[r]^{R}\ar@/^1cm/[rr]^{I}&Z^{p}(DD^{\per}(W))\ar[d]^{\int_{W/B}\hA(\hat o)\cup \dots}&\KU^{p}(W)\ar[d]^{\pi_{!}^{o}}\\\KU\C/\Z^{p-d-1}(B)\ar@/_1cm/[rr]&DD^{\per}(B)^{p-d-1}/\im(d)\ar[r]^(0.6){a}&\widehat{KU}^{p-d}(B)\ar[r]\ar@/_1cm/[rr]_{I}&Z^{p-d}(DD^{\per}(B))&\KU^{p-d}(B)}
\end{equation}

The set of $\widehat{KU}$-orientations refining an underlying topological $\KU$-orientation $o$ is a torsor over $DD^{\per}(W)^{-1}/\im(d)$ such that for $\alpha\in DD^{\per}(W)^{-1}/\im(d)$ we have
\begin{equation}\label{lwefwwfewfewfefewfewfewffwe}
\hA(\hat o+\alpha)=\hA(\hat o)+d\alpha .\end{equation}

In the case that $B$ is a point and $p-d$ is odd we have the rules (all following from \cite[(17)]{MR2664467})
\begin{equation}\label{kdklqwjdlqwjdwldqwdqwdqd}\pi_{!}^{\hat o+\alpha}(x)=\pi_{!}^{\hat o}(x)+\left[\int_{W}d\alpha\wedge R(x)\right]_{\C/\Z}\ , \quad  \pi_{!}^{\hat o}(a(\omega))=\left[\int_{W} \hA(\hat o)\wedge \omega\right]_{\C/\Z} \ ,\end{equation}
where we identify $\widehat{KU}^{p-d}(*)$ with $\C/\Z$, see \eqref{djhqkjwdqwdq9879}.

If $\pi=\pi_{1}\circ \pi_{0}$ is a composition of proper submersions and $\hat o_{i}$ are $\widehat{KU}$-orientations of $\pi_{i}$, then we can define a composed orientation $\hat o=\hat o_{1}\circ \hat o_{0}$ for $\pi$ in a natural way
such that
\begin{equation}\label{e23e23e3e23e32e32e2}
\pi_{!}^{\hat o}=\pi^{\hat o_{1}}_{1,!}\circ \pi_{0,!}^{\hat o_{0}}\ .   
\end{equation}
If
$$\xymatrix{W^{\prime}\ar[d]^{\pi^{\prime}}\ar[r]^{g}&W\ar[d]^{\pi}\\B^{\prime}\ar[r]^{f}&B}$$
is a cartesian diagram and $\hat o$ is a differential orientation of $\pi$, then we can define
a $\widehat{KU}$-orientation $\hat o^{\prime}$ of $\pi^{\prime}$ such that
$$\pi^{\prime,\hat o^{\prime}}_{!}\circ g^{*}=f^{*}\circ \pi_{!}^{\hat o}\ .$$

In order to avoid the additional complexity of the choice of $\mathrm{Spin}^{c}$-extensions of connections on real vector bundles with $\mathrm{Spin}^{c}$-structures in the present paper we will work with $\mathrm{Spin}$-structures. 
If $T^{v}\pi$ has a spin structure, then it has an induced $\mathrm{Spin}^{c}$-structure, and a connection
$ \nabla^{T^{v}\pi}$ has a canonical $\mathrm{Spin}^{c}$-extension, which we take from now on.
If the $\widehat{KU}$-orientation $\hat o$ is defined using the vertical metric and the horizontal distribution as above, then we have \begin{equation}\label{fwefw4343433fewfwf}
\hA(\hat o)=\hA(\nabla^{T^{v}\pi})\ ,\end{equation}
where $\nabla^{T^{v}\pi}$ is the Levi-Civita connection.

\bigskip

Assume that $\pi:W\to B$ is a submersion with fibrewise boundary $\partial \pi:\partial W\to B$. 
If $\hat o$ is a   $\widehat{KU}$-orientation of $\pi$, then we can define an induced $\widehat{KU}$-orientation $\partial \hat o$ of $\partial \pi$. In this situation we have the bordism formula \cite[Prop. 5.18]{MR2664467}. If $\hat x\in \widehat{KU}^{p}(W)$, then we have the equality \begin{equation}\label{jkhdkqjwdhqwdwqdqwdqwdqwd789}
\partial\pi^{\partial \hat o}_{!}(\hat  x_{|\partial W})=a\left(\int_{W/B} \hA(\hat o)\wedge R(\hat x )\right)\ .\end{equation}

\section{The invariant}\label{ergregregerger345345}

\subsection{Construction}\label{flkwjefwkejfeflkewjflkefjelkfjelfjewlfewf9790}

Given a closed odd-dimensional real foliated spin manifold $(M,\cF)$ such that $$2\codim(\cF)< \dim(M)$$ with a stable  framing $s$ of $\cF_{\R}$ and a complex vector bundle $(V,\nabla^{I})$  with flat partial connection we define an invariant $$\rho(M,\cF,\nabla^{I},s)\in \C/\Z\ .$$ 

\bigskip

In order to define the invariant we first choose the following additional geometric data:
\begin{enumerate}
\item We choose a connection $\nabla$ on $V$ which extends $\nabla^{I}$, see Definition \ref{fjwelfwfewf42343242wf}.
\item We choose an extension $\nabla^{\cF_{\R}^{\perp}}$ of the flat partial connection $\nabla^{I,\cF^{\perp}_{\R}}$, see Example \ref{kjlwefwefwef}.

\item We choose a Riemannian metric $g^{TM}$. 
\end{enumerate}

By Definition \ref{hffhweuifhiuwefwfewfwefew} a stable   framing $s$ of $\cF_{\R}$ is an isomorphism of real vector bundles
$$s:\cF_{\R}\oplus \underline{\R^{n}}\cong \underline{\R^{m}}$$
for certain choices of $n,m\in \nat$. The  trivial flat connection $\nabla^{s}$ on $\cF_{\R}\oplus \underline{\R^{n}}$ associated to the stable framing is induced from the trivial connection on $\underline{\R^{m}}$ via the isomorphism $s$.
The Riemannian metric on $M$ further induces an orthogonal splitting $TM\cong \cF_{\R}\oplus \cF_{\R}^{\perp}$ so that we can consider both connections
$\nabla^{LC}\oplus \nabla^{\underline{\R^{n}}}$ and $\nabla^{s}\oplus \nabla^{\cF^{\perp}_{\R}}$ on the same bundle $$TM\oplus \underline{\R^{n}}\cong (\cF_{\R}\oplus \underline{\R^{n}})\oplus \cF_{\R}^{\perp}\ .$$ In particular, we can define the transgression
 form (see  \eqref{nkdqlwdwqddwqdj09809}) 
$$\widetilde{\hA}(LC,s):=\widetilde{\hA}(\nabla^{LC}\oplus \nabla^{\underline{\R^{n}}},\nabla^{s}\oplus \nabla^{\cF^{\perp}_{\R}})\in DD^{\per}(M)^{-1}/\im(d)\ ,$$
where $\nabla^{LC}$ is the Levi-Civita connection on $TM$ associated to the Riemannian metric $g^{TM}$.

\bigskip

We now consider the map $\pi:M\to *$. Since $M$ is closed this is a proper submersion.
Since $M$ is spin, this map has a $KU$-orientation $o$. The choice of a Riemannian metric
refines the orientation $o$ to a  $\widehat{KU}$--orientation $\hat o$ (note that
the horizontal bundle is the zero bundle), see Subsection \ref{wlekfjwelfjewlfwef123}.  
 
 \begin{ddd}\label{flwefjwefewff} Let
  $M$ be an odd-dimensional real foliated closed spin manifold, let $s$ be a stable framing of $\cF_{\R}$, and let $\nabla^{I}$ be a flat partial connection on a complex vector bundle on $M$. Assume further that we have fixed $g^{TM}$, $\nabla^{\cF_{\R}^{\perp}}$,  and $\nabla$.
 Then
 we define
 $$\rho(M,\cF,\nabla^{I},s):=\pi_{!}^{\hat o-\widetilde{\hA}(LC,s)}([V,\nabla])\in \widehat{KU}^{-\dim(M)}(*)\stackrel{\eqref{djhqkjwdqwdq9879}}{\cong} \C/\Z\ .$$
\end{ddd}

In general, this quantity depends on the additional choices $g^{TM}$, $\nabla^{\cF_{\R}^{\perp}}$,  and $\nabla$. It will be a consequence of the bordism invariance that $\rho(M,\cF,\nabla^{I},s)$ is actually independent of these choices provided $2\codim(\cF)<\dim(M)$.

\begin{prop} \label{fklwfjeklfwefewfwf897}
Assume that $(M,\cF)$, $\nabla^{I}$, $s$ as well as  $g^{TM}$, $\nabla^{\cF_{\R}^{\perp}}$,  and $\nabla$ are as in Definition \ref{flwefjwefewff} with the   exception that $M$ is even-dimensional and has a boundary $\partial M$ which is transversal to $\cF$. We further assume that the geometric structures have a product structure near $\partial M$.
Then we have the equality
  $$\rho(\partial M,\cF_{|\partial M},\nabla^{I}_{ |\partial M},s_{|\partial M})=\left[\int_{M}\hA^{-}(\nabla^{\cF_{\R}^{\perp}})\wedge \ch^{-}(\nabla)\right]_{\C/\Z}\ .$$
  In particular, if $2\codim(\cF)<\dim(M)$, then $\rho(\partial M,\cF_{|\partial M},\nabla^{I}_{ |\partial M},s_{|\partial M})=0$.
\end{prop}
\proof
By the bordism formula \eqref{jkhdkqjwdhqwdwqdqwdqwdqwd789} we have
$$\rho(\partial M,\cF_{|\partial M},\nabla^{I}_{ |\partial M},s_{|\partial M})=\left[\int_{M} \hA(\hat o-\widetilde{\hA}(LC,s))\wedge R([V,\nabla])\right]_{\C/\Z}\ .$$
Using \eqref{lwefwwfewfewfefewfewfewffwe}, \eqref{fwefw4343433fewfwf} and \eqref{dkdqwdwqdqwdwqdkl8997}
we get
$$\rho(\partial M,\cF_{|\partial M},\nabla^{I}_{ |\partial M},s_{|\partial M})=\left[\int_{M} (\hA(\nabla^{LC})-d\widetilde{\hA}(LC,s))\wedge \ch(\nabla) \right]_{\C/\Z}\ .$$
We apply \eqref{1e1h2ekj12ej12ke21hek2eh} and the fact that $\hA$ is multiplicative in order to rewrite this as
$$\rho(\partial M,\cF_{|\partial M},\nabla^{I}_{ |\partial M},s_{|\partial M})=\left[\int_{M} \hA(\nabla^{\cF_{\R}^{\perp}})\wedge \ch(\nabla)\right]_{\C/\Z}\ .$$
We now use that both, $\nabla^{\cF_{\R}^{\perp}}$ and $\nabla$ extend flat partial connections. The associated characteristic forms   therefore refine to cycles in $DD^{-}(M)$.
Hence by Example \ref{kjlwefwefwef} and Definition \ref{qldjqwldqwdqwdqwd} we have
$$ \int_{M}\hA(\nabla^{\cF_{\R}^{\perp}})\wedge \ch(\nabla)=\int_{M}\hA^{-}(\nabla^{\cF_{\R}^{\perp}})\wedge \ch^{-}(\nabla) \ .$$
This implies the first assertion.

\bigskip

The integral of $\int_{M}$ factorizes over the  component in $$DD^{-}(M)(p)^{0}=F^{p}\Omega^{2p}(M)$$ with $p=\dim(M)/2$. 
If $\codim(\cF)<p$, then  we have $F^{p}\Omega^{2p}(M,\cF)=0$ and hence \begin{equation}\label{rferfefref32234324344}
Z^{0}(DD^{-}(M)(2p))=0\ .
\end{equation}
This implies the second  claim.
\hB

In the following we define the opposite of a framing and a  spin structure. 
Let $(M,\cF)$ and the stable framing  $s$ of $\cF_{\R}$ be given. Then we can form the cylinder
$I\times M$ with the foliation $T_{\C}I\boxplus  \cF$, see Example \ref{dekldlqwdqwdqwd}. We trivialize $TI\cong I\times \underline{\R}$ using the section $\partial_{t}$, where $t$ is the standard coordinate of the cylinder.  Then we write
$T(I\times M)\cong TI\boxplus TM\cong \underline{\R}\boxplus TM$ in order to define the induced spin structure on $I\times M$. Furthermore, the identification $TI\boxplus  \cF_{\R}\cong \underline{\R}\boxplus \cF_{\R}$ provides the stable framing $I\times s$ of $TI\boxplus  \cF_{\R}$. These constructions are made such that
$(M,\cF,s)$ is the boundary of $(I\times M,T_{\C}I\boxplus  \cF,I\times s)$ at the upper face of the cylinder corresponding to $1\in I$. \begin{ddd} We define $(M^{\op},\cF,s^{\op})$ to be the boundary of the cylinder at $0\in I$. \end{ddd} Here $M^{\op}$ indicates that $M$ is equipped with the opposite spin structure.

\bigskip

We adopt all the assumptions made in Definition \ref{flwefjwefewff} and fix choices
for $\nabla$, $\nabla^{\cF_{\R}^{\perp}}$ and $g^{TM}$. These can be extended constantly over the cylinder.  In this case $\hA^{-}(\nabla^{I\times \cF_{\R}^{\perp}})\wedge \ch^{-}(\pr^{*}\nabla)$
is pulled back from $M$ and has no $dt$-component. Hence its integral over $I\times M$  vanishes.
\begin{kor} 
 $$\rho(M,\cF,\nabla^{I},s)=-\rho(M^{\op},\cF,\nabla^{I},s^{\op})\ .$$
\end{kor}

Assume now that we have two choices for $\nabla$, $\nabla^{\cF_{\R}^{\perp}}$ and $g^{TM}$.   Then we can again consider the cylinder over $M$ and interpolate between these choices. The second assertion of Proposition \ref{fklwfjeklfwefewfwf897}  and the vanishing \eqref{rferfefref32234324344} for $p=\frac{\dim(M)+1}{2}$ imply:
\begin{kor}\label{djwqdlqwjdkwqdwqdqd} If $2\codim(\cF)< \dim(M)$, then  
 $\rho(M,\cF,\nabla^{I},s)$ is independent of the choices of  $\nabla$, $\nabla^{\cF_{\R}^{\perp}}$ and $g^{TM}$.
\end{kor}

\subsection{A spectral geometric interpretation}\label{fkweflwefewfewfewfwfwfw}

In this section we express $\rho(M,\cF,\nabla^{I},s)$ in terms of  spectral invariants of Dirac operators.

\bigskip

Let $M$ be a closed spin manifold with a Riemannian metric $g^{TM}$ and let $\bV:=(V,h^{V},\nabla^{u})$ be  a hermitean vector bundle with a unitary connection.  Then we can form the twisted Dirac operator $\Dirac\otimes \bV$. It is a first order elliptic differential operator which acts on the space of sections of $\Gamma(M,S(TM)\otimes V)$, where $S(TM)$ is the spinor bundle of $M$. It is symmetric with respect to the natural $L^{2}$-metric. Its spectrum is real and consists of eigenvalues
of finite multiplicity accumulating at $\pm\infty$.  By Weyl's asymptotics the number of eigenvalues with absolute value $\le R$ (counted with multiplicity) grows as $R^{\dim(M)}$. 
The $\eta$-function of this operator was introduced by Atiyah-Patodi-Singer \cite{MR0397797} and is defined by
$$\eta(\Dirac\otimes \bV)(s)=\sum_{\lambda\not=0} m_{\lambda} \sign(\lambda)|\lambda|^{-s}\ ,$$
where the sum is taken over the non-zero eigenvalues of $\Dirac\otimes \bV$ and $m_{\lambda}$ denotes the multiplicity. The sum converges for $\Ree(s)>\dim(M)$.  It has been further shown by  \cite{MR0397797} that the $\eta$-function  has a meromorphic continuation to all of $\C$ which is regular at $s=0$. 

\begin{ddd} The $\eta$-invariant of $\Dirac\otimes \bV$ is defined by $$\eta(\Dirac\otimes \bV):=\eta(\Dirac\otimes \bV)(0)\ .$$ We further define the reduced $\eta$-invariant
$$\xi(\Dirac\otimes\bV):=\left[\frac{\eta(\Dirac\otimes \bV)+\dim(\ker(\Dirac\otimes \bV))}{2}\right]_{\C/\Z}\in \C/\Z\ .$$
\end{ddd} 

We consider the projection $\pi:M\to *$. The spin structure on $M$ and the Riemannian metric $g^{TM}$ induce a $\widehat{KU}$-orientation $\hat o$, see Subsection \ref{wlekfjwelfjewlfwef123}. The geometric bundle
$\bV$ defines a class $[V,\nabla^{u}]\in \widehat{KU}^{0}(M)$. Using the identification \eqref{djhqkjwdqwdq9879} we get by \cite[Cor. 5.5]{MR2664467}

\begin{prop}\label{kqwdkjqdwqdwd}
$$\pi_{!}^{\hat o}([V,\nabla^{u}])=\xi(\Dirac\otimes\bV)\ .$$
\end{prop}

We now adopt the assumptions of Definition \ref{flwefjwefewff}.
We further choose a hermitean metric $h^{V}$ and a unitary connection $\nabla^{u}$ and set $\bV:=(V,h^{V},\nabla^{u})$.

\begin{prop}\label{fjwfkljfklfjwelkfjkewjfewlkf9798237982749237432243}
We have the equality
$$\rho(M,\cF,\nabla^{I},s)=\xi(\Dirac\otimes\bV)-\left[\int_{M} \widetilde{\hA}(LC,s)\wedge \ch(\nabla)-
\int_{M}\hA(\nabla^{LC})\wedge \widetilde{\ch}(\nabla,\nabla^{u})\right]_{\C/\Z}\ .$$
\end{prop}
\proof
We use the rules \eqref{kdklqwjdlqwjdwldqwdqwdqd} for $x=[V,\nabla]$, $\alpha=-\widetilde{\hA}(LC,s)$ and $\omega=\widetilde{\ch}(\nabla,\nabla^{u})$, Proposition \ref{kqwdkjqdwqdwd},  and 
equations \eqref{dkdqwdwqdqwdwqdkl8997}, \eqref{fwefw4343433fewfwf} and \eqref{kdqklwdjqwlkdjqlwdkwqd89798}. \hB

\section{Special cases}\label{jbjkhwekfewfwefewf89798}

\subsection{Adams $e$-invariant}\label{kdljlqwdqwdqwd}

We consider the case of the maximal foliation $\cF_{\max}=T_{\C}M$ on a closed manifold $M$, see Example \ref{lkjwqdwdqwdqwdqwd}. Then a stable framing $s$ of $\cF_{max,\R}$  is a stable  framing (see Definition \ref{hffhweuifhiuwefwfewfwefew}) of $TM$ and
$(M,s)$ defines a framed bordism class $[M,s]\in \Omega^{\fr}_{\dim(M)}$. The Pontrjagin-Thom construction identifies the framed bordism     theory $\Omega^{\fr}_{*}$ with the  homology theory represented by the sphere spectrum. In particular, its coefficients are the stable homotopy groups of the sphere $\Omega^{\fr}_{*}(*)\cong \pi_{*}^{s}$. In his study of the j-homomorphism Adams defined  in \cite{MR0198470}  a homomorphism
$$e^{\Adams}_{\C}:  \pi^{s}_{k}\to \C/\Z$$
for odd $k\in \nat$. A spectral geometric interpretation of $e^{\Adams}_{\C}$ has been given by Atiyah-Patodi-Singer in  
\cite{MR0397797}. In the following we describe $e^{\Adams}_{\C}$ using differential $KU$-theory.

The stable framing $s$ induces a spin structure on $M$. Given a Riemannian metric
$g^{TM}$ we obtain a $\widehat{KU}$-orientation $\hat o$ of $\pi:M\to *$. It has been observed in \cite[Prop. 5.22]{MR2664467}) that the $\widehat{KU}$-orientation 
$$\hat o_{s}:=\hat o-\widetilde{\hA}(LC,s)$$ of $\pi$
does not depend on the choice of the Riemannian metric. 

Let $\beins\in \widehat{KU}(M)$ be the class of the trivial one-dimensional bundle $(\underline{\C},\nabla^{\triv})$. 
Then by \cite[Lemma. 5.24]{MR2664467}  we have using \eqref{djhqkjwdqwdq9879}

 \begin{equation}\label{rewrewrrwr324324234234}
e^{\Adams}_{\C}([M,s])= \pi_{!}^{\hat o_{s}}(\beins)\ .
\end{equation}
The following corollary immediately follows from the Definition \ref{flwefjwefewff}.
\begin{kor}
We have
$$\rho(M,\cF_{\max},\nabla^{\triv},s)= e^{\Adams}_{\C}([M,s])\ .$$
\end{kor}

 Using the first identity in \eqref{kdklqwjdlqwjdwldqwdqwdqd} we get
the expression (to be used later) \begin{equation}\label{kljklqwdwqd}e^{\Adams}_{\C}([M,s])= \pi_{!}^{\hat o}(\beins) -\left[\int_{M} \widetilde{\hA}(LC,s)\right]_{\C/\Z} \end{equation} for the $e$-invariant.

\subsection{The $\rho$-invariant of flat bundles}
  
Assume now that $(V,\nabla)$ is a flat bundle on a  closed odd-dimensional spin manifold $M$. The spin structure on $M$ equips the map $\pi:M\to *$ with  an   orientation $o$ for $\KU$. We   observe that $$[V,\nabla ]-\dim(V)\beins\in \widehat{KU}^{0}_{\flat}(M)\cong \bKU\C/\Z^{-1}(M)\ .$$
 Hence we can apply the  integration map  \eqref{hjkfhjkhkjhkjefwefefwe7987} to this difference.
\begin{ddd}
We define the $\rho$-invariant of $\nabla$ by 
 $$\rho(\nabla):=\pi_{!}^{o}([V,\nabla]-\dim(V)\beins)\in \C/\Z \ .$$
\end{ddd}

If we choose a Riemannian metric, then we get a refinement of $o$ to a $\widehat{KU}$-orientation $\hat o$ of $\pi$. Using the integration in differential $K$-theory  and \eqref{r23r23r23r23r23r235435346546} we can write
\begin{equation}\label{jhjkhjkehfkewhfkfhkewjhfkewjfewfwf}
\rho(\nabla):=\pi_{!}^{\hat o}([V,\nabla])-\dim(V)\pi_{!}^{\hat o}(\beins)  \ .
\end{equation}

\begin{rem}{\rm 
Assume that $h^{V}$ is a hermitean metric on $V$ preserved by $\nabla$. Then we can form the geometric bundle $\bV=(V,\nabla,h^{V})$. As a consequence of Proposition \ref{kqwdkjqdwqdwd}   we have \begin{equation}\label{jhdjkehdkwwedewd}
\pi_{!}^{\hat o}([V,\nabla])-\dim(V)\pi_{!}^{\hat o}(\beins)=\xi(\Dirac\otimes \bV)-\dim(V)\xi(\Dirac)\ .
\end{equation}
Combining this with \eqref{jhjkhjkehfkewhfkfhkewjhfkewjfewfwf} we get the statement of the 
  index theorem for flat vector bundles by Atiyah-Patodi-Singer \cite[Thm. 5.3]{MR0397799}:   
$$ \xi(\Dirac\otimes \bV)-\dim(V)\xi(\Dirac)=\rho(\nabla)\ .$$
Observe that this is really a non-trivial statement. The left-hand side of this equality is the analytic index of the flat bundle, and the right-hand side is the topological index since we have defined 
 the $\rho$-invariant using the topological integration in $\KU\C/\Z$-theory.
}
\end{rem}

Let us now assume that the   spin structure on $M$ is induced by  a stable framing $s$ of $TM$.  
\begin{lem}
We have
\begin{equation}\label{mkmxlkmlqwxqwx}\rho(M,\cF_{\max},\nabla,s)=\rho(\nabla)+\dim(V) e^{\Adams}_{\C}([M,s])\ .\end{equation}
\end{lem}
\proof
Since $\nabla$ is flat we have
\begin{equation}\label{dqjwqwddqwdqwdwqd}\pi_{!}^{\hat o-\widetilde{\hA}(LC,s)}([V,\nabla])=\pi_{!}^{\hat o}([V,\nabla])-\dim(V)\left[\int_{M} \widetilde{\hA}(LC,s)\right]_{\C}\ .\end{equation}
We first use \eqref{kljklqwdwqd}
 in order to replace the second term in \eqref{dqjwqwddqwdqwdwqd} and then apply \eqref{jhjkhjkehfkewhfkfhkewjhfkewjfewfwf}.  \hB 

\begin{rem}{\rm 
The decomposition \eqref{mkmxlkmlqwxqwx} of the invariant $\rho(M,\cF_{\max},\nabla,s)$ is very interesting. 
A priori this quantity depends on the isomorphism class of the flat bundle $(V,\nabla)$.
But we now observe that  $\rho(M,\cF_{\max},\nabla,s)$ is actually an invariant of the class $[V,\nabla]^{\alg}\in \bK(\C)^{0}(M)$   represented by $(V,\nabla)$.  
This fact has already been shown in \cite{westburyjones} as we will explain in the following. 
 
 Fix a base point $m\in M$, choose an identification $V_{m}\cong \C^{\dim(V)}$,   and let $\alpha:\pi_{1}(M,m)\to GL(\dim(V),\C)$ denote the holonomy representation associated   to the flat connection $\nabla$ on $V$.
Then the quantity $e(M,\alpha)\in \C/\Z$ introduced by \cite{westburyjones} (for $M$ a homology sphere) can be written in the form (compare \cite[Thm. A]{westburyjones})
$$e(M,\alpha)=\rho(\nabla)\ .$$
The number  $e(M,\alpha)\in \C/\Z$ only depends on the algebraic $K$-theory class of $M$
determined by   $\alpha$,  which in our notation is $[V,\nabla]^{\alg}\in \bK(\C)^{0}(M)$.
Since clearly  $\dim(V)$ is an invariant of $[V,\nabla]^{\alg}$ as well, the combination
$$\rho(M,\cF_{\max},\nabla,s)=e(M,\alpha)+\dim(V) e^{\Adams}_{\C}([M,s])$$
only depends on the class $[V,\nabla]^{\alg}$ of $(V,\nabla)$.
 
In  Section \ref{keklwfewfewfewf} we will show a much stronger result. We will see that $ \rho(M,\cF_{\max},\nabla,s)$
only depends on the class  $$\pi^{o_{s}}_{!}([V,\nabla]^{\alg})\in \bK(\C)^{-\dim(M)}(*)\cong  K_{\dim(M)}(\C)\ ,$$ where $\bK(\C)^{*}$ is the cohomology theory represented by the algebraic $K$-theory spectrum of $\C$, and $o_{s}$ is the orientation of $\pi:M\to *$   for stable cohomotopy (and  hence for every cohomology theory since it is a module theory over stable cohomotopy)  given by the framing $s$.


The formula \eqref{mkmxlkmlqwxqwx} can be compared with formulas in Theorem \cite[Thm 5.5]{2011arXiv1103.4217B}. We conclude that $\rho(M,\cF_{\max},\nabla,s)$ can be expressed in terms of the universal $\eta$-invariant introduced in that reference. }
\end{rem}

 \subsection{A families $e$-invariant}

Let $$q:W\to B$$ be a proper submersion of relative dimension $p:=\dim(W)-\dim(B)>0$ and consider the vertical foliation $\cF^{v}=T^{v}_{\C}q$ on $W$, see Example \ref{djqlwdqwdqwdq}. This foliation is real. A  framing $s$ of $\cF_{\R}^{v}$ induces an orientation $o_{q}$ of the map $q$ for the framed bordism cohomology theory $\Omega^{\fr,*}$. We get a class $$[W\stackrel{q}{\to} B,s]=q^{o_{q}}_{!}(1_{S})\in \Omega^{\fr,-p}(B)\ ,$$ where $1_{S}\in \Omega^{\fr,0}(W)$ is the unit. The construction \eqref{rewrewrrwr324324234234}   of Adams' $e$-invariant can be extended from $B=*$ to general $B$ 
as a map
$$e^{\Adams}_{\C}:\Omega^{\fr,-p}(B)\to \KU\C/\Z^{-p-1}(B)\ .$$
According to \cite[Definition 5.23]{MR2664467} its value on the class $[W\stackrel{q}{\to} B,s]$ is given by
 \begin{equation}\label{ewfwfewfwfwf23434}
 e^{\Adams}_{\C}([W\stackrel{q}{\to} B,s]):=q_{!}^{\hat o_{s}}(\beins)\in \widehat{KU}^{-p}(B)_{\flat}\cong \KU\C/Z^{-p-1}(B)\ , 
\end{equation}
where $\hat o_{s}$ is the $\widehat{KU}$-orientation of $q$ induced by the framing, see Remark \ref{jdlqwkdjlqwdwqd}.

\begin{rem}\label{jdlqwkdjlqwdwqd}{\rm
The construction of $\hat o_{s}$ is similar as in Subsection \ref{kdljlqwdqwdqwd}. 
The vertical framing induces a spin structure. We choose a fibrewise Riemannian metric and a horizontal distribution.
Then we get a vertical Levi-Civita connection $\nabla^{T^{v}q}$. As explained in 
Subsection \ref{wlekfjwelfjewlfwef123} we get a $\widehat{KU}$-orientation $\hat o$. Furthermore, using the trivial
connection induced by the framing, we can define the transgression
$$\widetilde{\hA}(\nabla^{T^{v}q},s)\in DD^{\per}(W)^{-1}/\im(d)\ .$$
The $\widehat{KU}$-orientation
$$\hat o_{s}:=\hat o-\widetilde{\hA}(\nabla^{T^{v}q},s)$$ is then independent of the choice of the geometric structures.

In order to see that $e^{\Adams}_{\C}([W\stackrel{q}{\to} B,s])$ is flat we calculate its curvature using \eqref{r23r23r23r23r23r235435346546} and \eqref{lwefwwfewfewfefewfewfewffwe}
$$R(e^{\Adams}_{\C}([W\stackrel{q}{\to} B,s]))=\int_{W/B}(\hA(\hat o_{s})-d  \widetilde{\hA}(\nabla^{T^{v}q},s))=0\ .$$

}\end{rem}

Let us now assume that $B$ is closed and has a spin structure. Then the projection $\pi:B\to *$ has a $KU$ orientation $o_{\pi}$.
We choose a Riemannian metric $g^{TB}$ on $B$, a vertical metric $g^{T^{v}q}$, and a horizontal distribution $T^{h}q$. The metric $g^{TB}$ lifts to a metric on the horizontal bundle $T^{h}q$ and induces, together with the vertical metric $g^{T^{v}q}$,  a metric on $W$. Furthermore, the spin structure of $B$ induces a spin structure on the horizontal bundle, which together with the framing of $T^{v}q$ provides a spin structure on $W$.  Note that $\cF^{\perp}_{\R}\cong T^{h}\pi$. The Levi-Civita connection of $g^{TB}$ pulls back to the connection $\nabla^{\cF^{\perp}_{\R}}$.

\bigskip

We consider a geometric vector bundle $(V,\nabla)$ on $B$.
Then $(\pi^{*}V,\pi^{*}\nabla)$ is a bundle on $W$ and the restriction of $\pi^{*}\nabla$ to $\cF$ is flat, see Example \ref{fjewflewjflewfewf}.
 
 \bigskip
 
We now assume that $\dim(W)$ is odd.
\begin{lem}
We have the equality
$$\rho(W,\cF^{v},\pi^{*}\nabla,s)=\pi^{o_{\pi}}_{!}(e^{\Adams}_{\C}([W\stackrel{q}{\to}B,s])\cup [V])\ .$$
\end{lem}
\proof
The geometry on $B$ provides a $\widehat{KU}$-orientation $\hat o_{\pi}$. The geometry on $W$ induces a $\widehat{KU}$-orientation $\hat o_{\pi\circ q}$. In the following calculation we use \cite[Definition 3.22]{MR2664467}   and $\hA(\hat o_{\pi})=\hA(\nabla^{\cF_{\R}^{\perp}})$ at the place marked by $!$.
\begin{eqnarray*} \hat o_{\pi}\circ \hat o_{s}&=&\hat o_{\pi}\circ (\hat o-\widetilde{\hA}(\nabla^{T^{v}q},s))\\&=&\hat o_{\pi}\circ \hat o-\hA(\hat o_{\pi})\wedge \widetilde{\hA}(\nabla^{T^{v}q},s)\\&\stackrel{!}{=}&\hat o_{\pi\circ q}-\hA(\nabla^{\cF_{\R}^{\perp}})\wedge \widetilde{\hA}(\nabla^{T^{v}q},s)-\widetilde{\hA}(\nabla^{LC} , \nabla^{T^{v}q}\oplus  \nabla^{\cF_{\R}^{\perp}}  )\\&=&\hat o_{\pi\circ q}-\widetilde{\hA}(\nabla^{T^{v}q}\oplus  \nabla^{\cF_{\R}^{\perp}},  \nabla^{\mathrm{triv},s}\oplus \nabla^{\cF^{\perp}_{\R}})-
\widetilde{\hA}(\nabla^{LC} , \nabla^{T^{v}q}\oplus  \nabla^{\cF_{\R}^{\perp}}  ) \\&=&\hat o_{\pi\circ q}-
 \widetilde{\hA}(LC,s)\ .\end{eqnarray*}
We now use that integration 
is compatible with the  identification
$\KU\C/\Z^{*-1}\cong \widehat{KU}^{*}_{\flat}$. 
 We get, using the projection formula for the integration in  $\widehat{KU}$-theory,
\begin{eqnarray*}
\pi^{o_{\pi}}_{!}(e^{\Adams}_{\C}([W\to B,s])\cup [V])&=&
\pi^{\hat o_{\pi}}_{!}(e^{\Adams}_{\C}([W\to B,s])\cup [V,\nabla])
\\&\stackrel{\eqref{ewfwfewfwfwf23434}}{=}&
\pi^{\hat o_{\pi}}_{!}(q_{!}^{\hat o_{s}}(\pi^{*}[V,\nabla]))\\&\stackrel{\eqref{e23e23e3e23e32e32e2}}{=}&(\pi\circ q)^{\hat o_{\pi}\circ \hat o_{s}}(\pi^{*}[V,\nabla])\\
&=&(\pi\circ q)_{!}^{\hat o_{\pi\circ q}-
 \widetilde{\hA}(LC,s)}(\pi^{*}[V,\nabla])\\
 &=&\rho(M,\cF,\pi^{*}\nabla,s)\ .
\end{eqnarray*} \hB

\subsection{The dependence on the framing}

Let $s,s^{\prime}$ be two stable framings of a foliation $\cF_{\R}$.  Then we get two connections
$\nabla^{s}$ and $\nabla^{s^{\prime}}$ on $\cF_{\R}\oplus \R^{n}$. Since these connections are flat, by  \eqref{1e1h2ekj12ej12ke21hek2eh} we get a cohomology class
$$\widetilde{\hA}(\nabla^{s^{\prime}},\nabla^{s})\in H^{-1}(DD^{\per}(M))\ .$$

\begin{ddd}
For every class $u\in \KU^{0}(M)$ we
define the relative $e$-invariant of the pair
$(s^{\prime},s)$ of stable framings of $\cF_{\R}$  by $$e_u(s^{\prime},s)=\left[\int_{M}\widetilde{\hA}(\nabla^{s^{\prime}},\nabla^{s})\cup \ch(u)\right]\in \C/\Z\ .$$
\end{ddd}

\begin{rem}{\rm 
If $\cF=\cF_{\max}$, then
$$e_{1}(s^{\prime},s)=e^{\Adams}_{\C}([M,s^{\prime}])-e^{\Adams}_{\C}([M,s])\ .$$
In this case $e_{1}(s^{\prime},s)$ takes values in the well-known finite subgroup $\im(e^{\Adams}_{\C})\subseteq \C/\Z$ calculated  by Adams.  
 }
\end{rem}

The proof of the following proposition is a straightforward calculation.
\begin{prop}
We adopt the assumptions of Definition \ref{flwefjwefewff} and assume that $s,s^{\prime}$ are stable framings of $\cF_{\R}$.
Then we have
$$\rho(M,\cF,\nabla^{I},s^{\prime})-\rho(M,\cF,\nabla^{I},s)=e_{[V]}(s^{\prime},s)\ .$$
\end{prop}

\subsection{Real and imaginary parts}

\subsubsection{The decomposition}
In this subsection we discuss the components $\rho(\dots)^{\R/\Z}$ and $\rho(\dots)^{i\R}$ of
$\rho(M,\cF,\nabla^{I},s)$  associated to the decomposition of the target group $$\C/\Z\cong \R/\Z\oplus i\R\ , \quad x=x^{\R/\Z}+x^{i\R}$$ into the real and the imaginary parts.

\bigskip

We adopt the assumptions made in Definition \ref{flwefjwefewff}.  In addition we choose a hermitean metric $h^{V}$ on the complex vector bundle $V$. Then we can define the adjoint  connection  $\nabla^{*}$ of $\nabla$ (see Remark \ref{dkjqwlkdjqwlkjdlwqkdwdwqdqwd}) and its unitarization   $$\nabla^{u}:=\frac{1}{2}(\nabla+\nabla^{*})$$  with respect to $h^{V}$.
   
  \bigskip
  
We use \eqref{kdqklwdjqwlkdjqlwdkwqd89798} in order to write \begin{equation}\label{gdhqgdhjgwqdjhgqwd7987}
[V,\nabla]=[V,\nabla^{u}]+a(\widetilde{\ch}(\nabla,\nabla^{u}))\ .
\end{equation}  
Then we calculate
\begin{eqnarray}\lefteqn{\rho(M,\cF,\nabla^{I},s)}&&\nonumber\\&\stackrel{\eqref{kdklqwjdlqwjdwldqwdqwdqd}}{=}&\pi^{\hat o}_{!}([V,\nabla])-\left[\int_{M} \widetilde{\hA}(LC,s)\wedge \ch(\nabla)\right]_{\C/\Z}\nonumber\\&\stackrel{\eqref{gdhqgdhjgwqdjhgqwd7987}}{=}&\pi^{\hat o}_{!}([V,\nabla^{u}])+\left[\int_{M} \hA(\nabla^{LC})\wedge \widetilde{\ch}(\nabla,\nabla^{u})\right]_{\C/\Z} -   \left[\int_{M}   \widetilde{\hA}(LC,s)\wedge \ch(\nabla)\right]_{\C/\Z}\nonumber\\
&\stackrel{\eqref{hfjwkjehfkjwehfewfewf897987}}{=}&\pi^{\hat o}_{!}([V,\nabla^{u}])+\left[\int_{M}  \hA(\nabla^{LC})\wedge \widetilde{\ch}(\nabla,\nabla^{u})\right]_{\C/\Z}\nonumber\\&& -  \left[  \int_{M}      \widetilde{\hA}(LC,s)\wedge (\ch(\nabla^{u})+d\widetilde{\ch}(\nabla,\nabla^{u})) \right]_{\C/\Z}\nonumber\\&=&
\pi^{\hat o}_{!}([V,\nabla^{u}])-\left[\int_{M}   \widetilde{\hA}(LC,s)\wedge \ch(\nabla^{u})  \right]_{\C/\Z} +\left[\int_{M} \hA(\nabla^{\cF^{\perp}_{\R}})\wedge \widetilde{\ch}(\nabla,\nabla^{u})\right]_{\C/\Z} \label{dlkqwdqwdqwdqwdwqdqw}\end{eqnarray}
using partial integration and Stokes' theorem in the last step.
 
The first two summands in  \eqref{dlkqwdqwdqwdqwdwqdqw} are real.  The following is the decomposition of the transgression Chern form  into the real and imaginary part (we use \eqref{fewfwefwefwefew32434234234} and \eqref{qwdqwdqwdwqdqwdwqdqwdqwd}):
$$ \widetilde{\ch}(\nabla,\nabla^{u})=\frac{ \widetilde{\ch}(\nabla,\nabla^{u})+ \widetilde{\ch}(\nabla^{*},\nabla^{u})}{2}+\frac{\widetilde{\ch}(\nabla,\nabla^{*})}{2}\ .$$
We get 
\begin{eqnarray}
\rho(M,\cF,\nabla^{I},s)^{\R/\Z}&=&\pi^{\hat o}_{!}([V,\nabla^{u}])-\left[\int_{M}   \widetilde{\hA}(LC,s)\wedge \ch(\nabla^{u})  \right]_{\R/\Z} \\&&+\left[\int_{M} \hA(\nabla^{\cF^{\perp}_{\R}})\wedge \frac{ \widetilde{\ch}(\nabla,\nabla^{u})+ \widetilde{\ch}(\nabla^{*},\nabla^{u})}{2}\right]_{\R/\Z}\nonumber\\[1cm]
 \rho(M,\cF,\nabla^{I},s)^{i\R}  &=&  \int_{M} \hA(\nabla^{\cF^{\perp}_{\R}})\wedge \frac{\widetilde{\ch}(\nabla,\nabla^{*})}{2} \label{fjilfewfwefffef}
\ .\end{eqnarray}

\subsubsection{The imaginary part}

We see that the imaginary part  $\rho(M,\cF,\nabla^{I},s)^{i\R}$ is just a characteristic number which can be calculated as an integral over locally computable quantities. It does not depend on the framing.

\begin{ex}{\rm  We assume that $\nabla^{I}$ is unitary with respect to the metric $h$. Then we can take for $\nabla$ the unitary extension constructed in Lemma 
\ref{dhqwkdqkwddqwdwqdioipopioipoopi}. With this choice we have $\nabla= \nabla^{u}$.  
\begin{kor}
If $\nabla$ is  the unitary extension of $\nabla^{I}$, then $$\rho(M,\cF,\nabla^{I},s)^{i\R}=0\ .$$
In particular, if $2\codim(\cF)<\dim(M)$ and $\nabla^{I}$ is unitary, then $\rho(M,\cF,\nabla^{I},s)^{i\R}=0$.  \end{kor}\proof
The first assertion  follows from \eqref{fjilfewfwefffef}, $\nabla=\nabla^{*}$ and the second equality in \eqref{fewfwefwefwefew32434234234}. The second assertion  is then a consequence of the first and Corollary \ref{djwqdlqwjdkwqdwqdqd}. \hB
}\end{ex}

\bigskip

 \begin{ex}{\rm If $\pi:\tilde M\to M$ is a finite covering of degree $[\tilde M:M]\in \nat$, then
we have the identity
 \begin{equation}\label{dwddwqqdqwd342423424324dasd}
\rho(\tilde M,\pi^{*}\cF,\pi^{*}\nabla^{I},\pi^{*}s)^{i\R}=[\tilde M:M] \rho(M,\cF,\nabla^{I},s)^{i\R}\ .
\end{equation}
}\end{ex}

\bigskip

Given a foliated manifold $(M,\cF)$ we have an associated bundle $\cF^{\perp}$ with a flat partial connection $\nabla^{I,\cF}$.  If we apply $\rho(\dots)^{i\R}$ to  $(V,\nabla^{I})=(\cF^{\perp},\nabla^{I,\cF})$ or a bundle obtained from this by some operation of tensor calculus we get an invariant of the foliation $(M,\cF)$.

\begin{ex}\label{fkllwefwefewf}{\rm

In this example, for even $n\in \Z$,  we consider a $2n+1$-dimensional closed oriented manifold $M$ with a real foliation $\cF$ of codimension $1$.  We assume that $\cF_{\R}$ is co-oriented. 

\bigskip

We first recall the definition of the Godbillon-Vey class $\mathbf{GV}_{2k+1}(\cF)\in H^{2k+1}(M;\R)$ for $k\ge 1$.
Since $\cF_{\R}$ is co-oriented there exists a real nowhere vanishing one-form $\kappa\in \Omega^{1}(M)$ such that 
$\cF_{\R}=\ker(\kappa)$. Integrability of $\cF_{\R}$ translates to the relation $\kappa\wedge d\kappa=0$. We can choose a real $1$-form
$\omega\in \Omega^{1}(M)$ such that $d\kappa=\kappa\wedge \omega$. 
Note that $ \omega$ is unique up to multiples of $\kappa$. Then the form $\omega\wedge (d\omega)^{k}\in \Omega^{2k+1}(M)$ is closed and represents the 
Godbillon-Vey class $\mathbf{GV}_{2k+1}(\cF)$.  

\bigskip

We now assume that  $TM$ has a stable framing $s_{M}$   and a Riemannian metric $g^{TM}$. 
The co-orientation of $\cF_{\R}$ induces a framing $s^{\perp}$ of $\cF^{\perp}_{\R}$ by the positive normal unit vector field $N$.
Then there is   a unique stable framing $s$ of $\cF_{\R}$ (up to equivalence and homotopy) such that $s\oplus s^{\perp}\sim s_{M}$.
For $(V,\nabla^{I})$ we take $(\cF^{\perp},\nabla^{I,\cF^{\perp}})$. 

\begin{lem}
Assume:
\begin{enumerate}
\item $M$ is a closed oriented $2n+1$-dimensional Riemanniann manifold with a stable framing $s_{M}$.
\item  $\cF$ is a real, co-oriented codimension-one foliation on $M$ with a stable framing $s$ of $\cF_{\R}$.
\item $(V,\nabla^{I})=(\cF^{\perp},\nabla^{I,\cF^{\perp}})$.
\item   We have 
$s\oplus s^{\perp}\sim s_{M}$,
where  $s^{\perp}$ is the stable framing of $\cF_{\R}^{\perp}$
given by the positive unit normal vector field.
\end{enumerate}

Then 
we have
$$\rho(M,\cF,\nabla^{I},s)^{i\R}=\frac{(-1)^{n+1}}{(2\pi i)^{n+1} n!}\int_{M} \mathbf{GV}_{2n+1}(\cF)\ .$$
\end{lem}
\proof Since $\dim(\cF_{\R}^{\perp})=1$ we
 have  $ \hA_{4p}(\nabla^{\cF^{\perp}_{\R}})\in F^{2p}\Omega^{4p}(M,\cF)=0$ for all $p\ge 1$.
  Hence   \eqref{fjilfewfwefffef} specializes to 
$$\rho(M,\cF,\nabla^{I,\cF^{\perp}},s)^{i\R}=\frac{1}{2}\int_{M} \widetilde{\ch}_{2n+2}(\nabla,\nabla^{*})\ .$$
So we must identify $\widetilde{\ch}_{2n+2}(\nabla,\nabla^{*})$ with a multiple of  $\mathbf{GV}_{2n+1}(\cF)$.

\bigskip

  Using the unit normal vector field $N\in \Gamma(M,TM)$ 
 we can normalize $\kappa$ such that $ \kappa(N)=1$. Let $\omega$ be as above. We take $\omega$ as a connection one-form for a connection $\nabla$ on $\cF^{\perp}_{\R}$ with respect to the trivialization by $N$.
For a section $X$ of $TM$ we have by definition
$$\nabla_{X}N=\omega(X)N\ .$$
On the other hand, if $X$ is a section of $\cF$, then we have by Cartan's formula
  $$\omega(X)=(\kappa\wedge \omega)(N,X)=d\kappa(N,X)= N\kappa(X)-X\kappa(N)-\kappa([N,X])=\kappa([X,N])\ .$$  
  In view of the description of $\nabla^{I,\cF^{\perp}_{\R}}$ given in Example \ref{fklwefjwefewf} this implies that  the connection $\nabla$ extends the flat partial connection $\nabla^{I,\cF^{\perp}_{\R}}$.   
  
  \bigskip
  
 We have
 $$\frac{(-1)^{n+1}}{(2\pi i)^{n+1} n!}\omega\wedge (d\omega)^{n}= \widetilde{\ch}_{2n+2}(\nabla,\nabla^{\triv})\ .$$
 Similarly, 
 $$(-1)^{n+1}\frac{(-1)^{n+1}}{(2\pi i)^{n+1} n!}\omega\wedge (d\omega)^{n}= \widetilde{\ch}_{2n+2}(\nabla^{*},\nabla^{\triv})\ .$$
 Hence, if $n$ is even, then by taking the difference of these two equations we get
 $$\frac{2(-1)^{n+1}}{(2\pi i)^{n+1} n!}\mathbf{GV}_{2n+1}(\cF)=  \widetilde{\ch}_{2n+2}(\nabla ,\nabla^{*})\ .$$
 
   \hB
 } 
 \end{ex}

 \begin{rem}\label{ergegojerglerogergeg}{\rm 
As noted  above we can take $(V,\nabla^{I}):=(\cF^{\perp},\nabla^{I, \cF^{\perp}})$
in order to define an invariant which only depends on the foliation $\cF$.
In this example assume that $\cF$ is real and that
 $\nabla^{\cF^\perp}$ is  the complexification of a connection $\nabla^{\cF_{\R}^{\perp}}$ extending $\nabla^{I,\cF^{\perp}_{\R}}.$ We choose in addition a metric $h^{\cF^{\perp}_{\R}}$ in order to define the adjoint $\nabla^{\cF^{\perp},*} $.
In this remark we explain the place of $$\rho(M,\cF,\nabla^{I,\cF^{\perp}},s)^{i\R}= \int_{M} \hA(\nabla^{\cF^{\perp}_{\R}})\wedge \frac{\widetilde{\ch}(\nabla^{\cF^{\perp}},\nabla^{\cF^{\perp},*})}{2}$$ in the classification of foliation invariants defined in terms of secondary characteristic classes of foliations.

\bigskip

We start with the classification of characteristic forms for foliations of codimension $q\in \nat$ \cite{MR0307250}, see also \cite{MR512428}. Let 
  $q^{\prime}\in \nat$ be the greatest odd integer $\le q$.
One defines the commutative graded algebra
$$WO_{q}:=\R[\tilde c_{1},\dots,\tilde c_{q^{\prime}}]\otimes \R[c_{1}, ,\dots,c_{q}]^{\le 2q}\ ,$$
where the degrees of the generators are given by $$|\tilde c_{j}|=2j-1\ , \quad \mbox{$j$ odd, $\quad \quad $ and }\quad \quad \quad |c_{j}|=2j$$
and the superscript $[-]^{\le 2q}$ indicates that we take only polynomials of degree less than $2q$.

On this ring we consider the differential $d$ given by
$$d\tilde c_{j}:=c_{j}\ ,\quad dc_{j}=0\ .$$ 
The cohomology $H^{*}(WO_{q})$ of this DGA classifies secondary characteristic classes for foliations of codimension $q$.
For a cohomology class $[U]\in H^{q}(WO_{q})$ we let $\Delta([U])\in H^{*}(M;\R)$ denote  the corresponding cohomology class.

\bigskip

In the following we describe $\Delta$ on the form level. Since  $\nabla^{\cF^{\perp}}$ and $\nabla^{\cF^{\perp},*}$ are  complexifications  of  connections which are   dual to each other on a real bundle we have
$$\ch_{2n}(\nabla^{\cF^{\perp},*})=(-1)^{n} \ch_{2n}(\nabla^{\cF^{\perp}})\ .$$
By  \eqref{hfjwkjehfkjwehfewfewf897987} 
we  get for odd $n$
$$d \frac{1}{2i^{n}}\widetilde{\ch}_{2n}(\nabla^{\cF^{\perp}},  \nabla^{\cF^{\perp},*})=\frac{1}{i^{n}}\ch_{2n}(\nabla^{\cF^{\perp}})\ .$$
Therefore the
 connection $\nabla^{\cF^{\perp}}$ together with a choice of a metric $h^{\cF^{\perp}_{\R}}$ induces a map of commutative differential graded algebras 
$$\Delta_{(\nabla^{\cF^{\perp}},h^{\cF^{\perp}_{\R}})}:WO_{q}\to \Omega(M)\ ,$$  by $$\Delta_{( \nabla^{\cF^{\perp}} ,h^{\cF^{\perp}_{\R}})}(\tilde c_{n}):=\frac{1}{2i^{n}}\widetilde{\ch}_{2n}(\nabla^{\cF^{\perp}},  \nabla^{\cF^{\perp},*})\ , \quad \Delta_{( \nabla^{\cF^{\perp}} ,h^{\cF^{\perp}_{\R}})} (c_{n}):=\frac{1}{i^{n}}\ch_{2n}(\nabla^{\cF^{\perp}})\ .$$

Then for $[U]\in H^{*}(WO_{q})$ the characteristic class $\Delta([U])\in H^{*}(M;\R)$ of the foliation $\cF$ is given by 
 \begin{equation}\label{qwdqwdwqdqwdqwdqwd2312343534tfwrvwfvwfwef}
\Delta([U]):=[\Delta_{( \nabla^{\cF^{\perp}} ,h^{\cF^{\perp}_{\R}})}(U)] \ .
\end{equation}

\bigskip

There is a universal polynomial
$ A(c_{1},\dots,c_{q})\in \R[c_{1},\dots,c_{q}]^{\le 2q}$  such that
$$\hA(\nabla^{\cF^{\perp}_{\R}})^{\le 2q}= A(\ch_{2}(\nabla^{\cF^{\perp}}),\dots,\ch_{2q}(\nabla^{\cF^{\perp}}))\ .$$
We consider \begin{equation}\label{fwefwefewfewfewfwefew53453453455}
U:=\left[\left(\sum_{j=1, \mathrm{odd}}^{q^{\prime}} (-1)^{\frac{j+1}{2}} \tilde c_{j}\right)  A(c_{1},\dots,c_{q}) \right]_{\dim(M)}\in WO^{\dim(M)}_{q}\ . 
\end{equation}

If $2q<\dim(M)$, then $U$ is a cycle. \begin{lem} Let $\cF$ be a real foliation of codimension $q$ such that $2q<\dim(M)$. Then the  class $$[U]\in H^{\dim(M)}(WO_{q})$$
is the universal class classifying the imaginary part of $\rho(M,\cF,\nabla^{I,\cF^{\perp}},s)$.
\end{lem}
\proof The relation  
 $$\rho(M,\cF,\nabla^{I,\cF^{\perp}},s)^{i\R} =i\langle \Delta([U]),[M]\rangle$$ follows immediately from \eqref{fwefwefewfewfewfwefew53453453455}, the definition \eqref{qwdqwdwqdqwdqwdqwd2312343534tfwrvwfvwfwef} of $\Delta([U])$ and \eqref{fjilfewfwefffef}.  
 \hB 
 
 Let us assume that $p$ is odd and $2p-1> q$. Then $d\tilde c_{p}=0$ and we have the cohomology class $[\tilde c_{p}]\in H^{2p-1}(WO_{q})$. 
 If the foliation $\cF$ is real, then the characteristic class
 \eqref{kfkwejwlkefjlwekfjewfewfopipoi234} is given by \begin{equation}\label{ewfwefewfwf432342344123}
 [c_{2p-1}(\nabla^{I,\cF^{\perp}})]=2i^{p}\Delta[\tilde c_{p}]\ .
\end{equation}

 }
\end{rem}

\subsubsection{The real part}

 The real part $\rho(M,\cF,\nabla^{I},s)^{\R/\Z}$ is more complicated and of global nature. 
 A good case to look at is discussed in Example \ref{kdljlqwdqwdqwd}.

\begin{ex}{\rm 
The following example shows that $\rho(M,\cF,\nabla^{I},s)$ is not an integral over $M$ of locally determined quantities.
We consider the  manifold  $M:=S^{1}$ with the maximal foliation $\cF_{\max}=T_{\C}S^{1}$. The framing $s$ of $TS^{1}$ is the bounding framing so that $[S^{1},s]=0$ in $\Omega_{1}^{\fr}$. Furthermore we  let   $\bV(r):=(V,h,\nabla(r))$
be the flat line bundle with holonomy $\exp(2\pi i r)$ for $r\in [0,1)$. Then we can apply \eqref{mkmxlkmlqwxqwx} and \eqref{jhdjkehdkwwedewd} and get
$$\rho(S^{1},\cF_{\max},\nabla(r),s)=\rho(\nabla)=\xi(\Dirac\otimes \bV(r))-\xi(\Dirac)\ .$$
In this case the reduced $\eta$-invariant can be calculated explicitly. The result is
$$\xi(\Dirac\otimes \bV(r))=[-r]_{\C/\Z}\ .$$
Hence we get
$$\rho(S^{1},\cF_{\max},\nabla(r),s) =[-r]_{\C/\Z}\ .$$
In particular, our invariant  depends non-trivially on $r$.
The data $(S^{1},\cF_{\max},\nabla(r),s)$  for different $r$ are locally isomorphic.

\bigskip

Note that in this example the analogue of \eqref{dwddwqqdqwd342423424324dasd} nevertheless holds true.
  }
\end{ex}

\section{Factorization over algebraic $K$-theory of smooth functions}\label{keklwfewfewfewf}
 
 Let $P$ be  a closed $p$-dimensional manifold and let $s$ be a stable framing of $TP$. For a manifold $X$
we consider a product of foliated manifolds \begin{equation}\label{gergergergregegrgerg3453435}
(M,\cF):=(P\times X,T_{\C}P\boxplus 0)=(P,\cF_{\max})\times (X,\cF_{\min})
\end{equation}  and a pair $(V,\nabla^{I})$ of a complex vector bundle and a flat partial connection on $(M,\cF)$. We will show that the data represents an algebraic $K$-theory 
 class $$f^{o_{s}}_{!}(
 [V,\nabla^{I}]^{\alg})\in K_{p}(C^{\infty}(X))$$ of the ring $C^{\infty}(X)$. 
  If we assume that $X$ is closed, spin and that $\dim(X)<p$, then
 our main result is   the equality $$\rho(M,\cF,\nabla^{I},s)=\pi_{!}^{o}(\reg_{X}(f^{o_{s}}_{!}([V,\nabla^{I}]^{\alg})))\ ,$$
where $$\reg_{X}:K_{p}(C^{\infty}(X))\to \ku\C/\Z^{-p-1}(X)$$ is the regulator and
the map $$\pi:X\to *$$ has the $\ku$-orientation $o$ from the spin structure of $X$.

 \subsection{Statement of the result}\label{}

For manifolds $X$ and $P$ we consider the foliated manifold \eqref{gergergergregegrgerg3453435}. 
From the point of view of foliation theory it is very trivial. The leaves of the foliation on the product    $M=P\times X$    are just the submanifolds $P\times \{x\}$ for all $x\in X$.    
 \bigskip
 
 We assume that $P$ is closed and that the tangent bundle $TP$ of $P$ is equipped with a  stable framing $s$.  The framing $s$ induces an orientation $o_{s}$ of the map $f:P\to *$ for stable   cohomotopy theory, the cohomology theory represented by the sphere spectrum $\bS$ (or equivalently, the framed bordism theory).
 Any spectrum $\bE$ is a module spectrum over $\bS$. Consequently $f$ has an induced 
   orientation for the cohomology theory $\bE^{*}$ which we denote by   the same symbol $o_{s}$. 
 We have an 
 Umkehr or integration  map between cohomology groups
   $$ f_{!}^{o_{s}}:\bE^{*}(P)\to \bE^{*-p}(*)\ ,$$
   where $p:=\dim(P)$.  
   We will apply this to the  cohomology theory $\bK(C^{\infty}(X))^{*}$ represented by the connective algebraic $K$-theory spectrum $\bK(C^{\infty}(X))$ of the ring of complex-valued smooth functions on the manifold $X$.

\bigskip

We start with the class $$[V,\nabla^{I}]^{\alg}\in \bK(C^{\infty}(X))^{0}(P) $$ (see Definition \ref{xgrergegg}  for a technical description) represented by a pair $(V,\nabla^{I})$ of 
 a complex vector  bundle  on the foliated manifold \eqref{gergergergregegrgerg3453435} and  a flat partial connection.
We can form the algebraic $K$-theory class   \begin{equation}\label{ghjqwdgjqwdwqdwqdqdqdqd}
f^{o_{s}}_{!}([V,\nabla^{I}]^{\alg})\in \bK(C^{\infty}(X))^{-\dim(P)}(*)=K_{\dim(P)}(C^{\infty}(X))\ .
\end{equation}


We  now  assume that $X$ is closed and spin. We further assume that $\dim(P)+\dim(X)$ is odd and that   $\dim(X)<\dim(P)$, or equivalently, $$2\codim(\cF)<\dim(M)\ .$$ Then by Corollary \ref{djwqdlqwjdkwqdwqdqd}  the invariant $\rho(M,\cF,\nabla^{I},s)\in \C/\Z$ is well-defined and independent of additional geometric choices. 

The main result of the present section shows that  $\rho(M,\cF,\nabla^{I},s)$ can be expressed in terms of the class
 \eqref{ghjqwdgjqwdwqdwqdqdqdqd}. 
In greater detail, for every $n\in \nat$ with $n> \dim(X)$ we will construct, using methods from differential cohomology theory, a natural  regulator
$$\reg_{X}:K_{n}(C^{\infty}(X))\to \ku\C/\Z^{-n-1}(X)\ ,$$
see Definition \ref{klfwefewfewfwf}.
Let $\pi:X\to *$ be the projection. The spin structure on $X$ induces an orientation $o$ for the periodic complex topological $K$-theory $\KU$, and hence for the $\KU$-modules    $\ku$ and $\ku\C/\Z$. 
We use the isomorphisms
$$\ku^{k}\cong \left\{\begin{array}{cc}\Z&k\in 2\nat\\ 0&\mathrm{else}\end{array}\right.  \qquad \mbox{and} \qquad \ku\C/\Z^{k}\cong \left\{\begin{array}{cc}\C/\Z&k\in 2\nat\\ 0&\mathrm{else}\end{array}\right.$$
in order to interpret elements in $\ku\C/\Z^{2*}(*)$ (e.g., the left-hand side of \eqref{wefwefefewfwfewfwewefwefwef}) as elements of $\C/\Z$.
\begin{theorem}\label{flkfefwefwefewfef} We have the relation
\begin{equation}\label{wefwefefewfwfewfwewefwefwef}\pi^{o}_{!}(\reg_{X}(f^{o_{s}}_{!}([V,\nabla^{I}]^{\alg})))=\rho(M,\cF,\nabla^{I},s)\ .\end{equation}
\end{theorem}
The proof of this theorem will be finished in Subsection \ref{fwklfwfewfewfewf}. 

\begin{rem}{\rm 
Every class $x\in K_{*}(C^{\infty}(X))$ can be presented  in the form \eqref{ghjqwdgjqwdwqdwqdqdqdqd} for suitable stably framed manifolds $P$ and pairs $(V,\nabla^{I})$. Indeed, the class $x$ can be thought of as being represented by  
a map $x:S^{n}\to BGL(C^{\infty}(X) )^{+}$, where we consider $GL(C^{\infty}(X) )$ as a discrete group and $+$ stands for Quillen's $+$-construction. Using the standard stable framing $s_{\mathrm{can}}$ of $S^{n}$ the  triple $(S^{n},x,s_{\mathrm{can}})$ represents a framed bordism class
$[S^{n},x,s_{\mathrm{can}}]\in \Omega^{\fr}_{n} (BGL(C^{\infty}(X) )^{+})$. Since the $+$-construction map
$$p:BGL(C^{\infty}(X) )\to BGL(C^{\infty}(X) )^{+}$$ induces an isomorphism in generalized homology theories there exists a unique class
$[P,y,s]\in \Omega^{\fr}_{n}(BGL(C^{\infty}(X) ))$ such that $p_{*}([P,y,s])=[S^{n},x,s_{\mathrm{can}}]$. Since $P$ is compact, there exists a factorization of $y$ as
$$P\stackrel{\tilde y}{\to} BGL(N,C^{\infty}(X) )\to BGL(C^{\infty}(X) )$$ for a suitable $N\in \nat$. The map $\tilde y$
classifies a pair $(V,\nabla^{I})$ over $P\times X$ of an $N$-dimensional complex vector  bundle with a flat partial connection in the $P$-direction. We then have
$$f_{!}^{o_{s}}([V,\nabla^{I}]^{\alg})=x\ .$$ \hB
}
\end{rem}

\subsection{Algebraic $K$-theory sheaves}
\label{klfwjlefewfewf}

 We consider the site $\Mf_{\C\textrm{-}\fol}$ of pairs $(M,\cF)$ of manifolds $M$ with a  foliation $\cF$ and foliated maps (see Section \ref{fewl453534535435} for definitions).  
 The topology is given by open coverings. We have a morphism of sites  \begin{equation}\label{ewfwfwefewfewfwfwefwfe}
\Mf_{\C\textrm{-}\fol}\to \Mf
\end{equation} which forgets the foliations.

In the following we work in the framework of $\infty$-, or more precisely, of $(\infty,1)$-categories developed by Joyal, Lurie and others \cite{HTT}, \cite{HA}. We refer to \cite[Sec. 2.1]{Bunke:2014aa}, \cite[Sec. 2]{2013arXiv1311.3188B}
and \cite[Sec.4]{2012arXiv1208.3961B} for an introduction to the language as we will use it here
and for  further references. We will not discuss the size issues. They can be solved  in the standard way 
for the examples used in the present paper.

\bigskip

 For a presentable $\infty$-category $\bC$ and a site $\bM$ we consider the category
 $$\PSh_{\bC}(\bM):=\Fun(\bM^{op},\bC)$$ of $\bC$-values presheaves and its full subcategory of sheaves $\Sh_{\bC}(\bM)$.
 They are related by an adjunction \begin{equation}\label{dqdqwdwqdwqdwqdwqdqd}
L:\PSh_{\bC}(\bM)\leftrightarrows\Sh_{\bC}(\bM):\mathrm{inclusion}\ ,
\end{equation}
 
 where $L$ is called the sheafification.  
 \bigskip
 
 We consider the $1$-category of categories $\Cat$ with its cartesian symmetric monoidal structure.  
 For the class $W$   of categorical equivalences  we form the symmetric monoidal $\infty$-category
 $\Cat[W^{-1}]$. By $\CAlg(\Cat[W^{-1}])$ we denote the category of commutative algebras  in $\Cat[W^{-1}]$. 
 \begin{rem}\label{dkwqnmqlkwdwqdwqdqd}{\rm A commutative monoid can be considered as a symmetric monoidal category with only unit morphisms.
It is an object of $\CAlg(\Cat)$ and therefore represents one in $\CAlg(\Cat[W^{-1}])$. A general symmetric monoidal category has non-identity associator and commutativity constraints and is therefore not a commutative algebra in $\Cat$. But it  naturally represents an object in $\CAlg(\Cat[W^{-1}])$.
 
 }
 \end{rem}
   The objects of $\PSh_{\CAlg(\Cat[W^{-1}])} (\bM)$ are called
 symmetric monoidal prestacks. Similarly,  objects in $\Sh_{\CAlg(\Cat[W^{-1}])}(\bM)$ are called symmetric monoidal  stacks.
   
   \bigskip

 We consider the following four symmetric monoidal stacks on $\Mf$ or $\Mf_{\C\textrm{-}\fol}$    of vector bundles with additional structures. The monoidal structure is always given by the direct sum.

\begin{enumerate}
\item For a manifold $M$
we let $\Vect (M )$ denote the category of    vector bundles $V\to M$. A  map $f:M \to M^{\prime}$   induces a functor $f^{*}:\Vect (M^{\prime} )\to \Vect (M )$. We get a stack $\Vect $ on the site $\Mf $ with respect to the topology of open coverings.  
We use the same symbol for its pull-back to the site $\Mf_{\C\textrm{-}\fol}$ along \eqref{ewfwfwefewfewfwfwefwfe}.
\item We let $\Vect^{\nabla} (M)$ denote the category of  pairs  ($V,\nabla)$ of a vector bundle  $V\to M$ and a connection. A   map $f: M  \to M^{\prime} $   induces a functor $f^{*}:\Vect^{\nabla} (M^{\prime} )\to \Vect^{\nabla} (M )$. We get a symmetric monoidal stack $\Vect^{\nabla} $ on  the site $\Mf $.  We use the same symbol for its pull-back to the site $\Mf_{\C\textrm{-}\fol}$ along \eqref{ewfwfwefewfewfwfwefwfe}.
\item For a foliated manifold $(M,\cF)$ we let $\Vect^{\flat}(M,\cF)$ denote the category of pairs $(V,\nabla^{I})$ of a vector bundles $V\to M$ and a flat partial connection $\nabla^{I}$ on $V$, see Section \ref{fhfjlwefkjfewfewfewfwf}. A foliated map     $f:(M,\cF)\to (M^{\prime},\cF^{\prime})$   induces a functor $f^{*}:\Vect^{\flat}(M^{\prime},\cF^{\prime})\to \Vect^{\flat}(M,\cF)$. We get a stack $\Vect^{\flat}$ on  the site $\Mf_{\C\textrm{-}\fol}$.  \item
We let $\Vect^{\flat,\nabla}(M,\cF)$ denote the category of pairs $(V,\nabla)$ of a vector bundle $V\to M$ and a   connection $\nabla$ on $V$ which is flat in the direction of the foliation. A foliated map $f$ as above induces a functor $f^{*}:\Vect^{\flat,\nabla}(M^{\prime},\cF^{\prime})\to \Vect^{\flat,\nabla}(M,\cF)$. We get a symmetric monoidal stack $\Vect^{\flat,\nabla}$ on  the site $\Mf_{\C\textrm{-}\fol}$.   \end{enumerate}

We will consider $\Vect$ and $\Vect^{\nabla}$ also as stacks on $\Mf_{\C\textrm{-}\fol}$  via pull-back along the forgetful morphism \eqref{ewfwfwefewfewfwfwefwfe}.
There is a commutative diagram of  forgetful  maps \begin{equation}\label{wsdqwdwqdqwdwqdsq}   \xymatrix{\Vect^{\flat,\nabla}\ar[d]\ar[r]&\Vect^{\flat}\ar[d]\\\Vect^{\nabla}\ar[r]&\Vect}\end{equation}  in  $\Sh_{\CAlg(\Cat[W^{-1}])} (\Mf_{\C\textrm{-}\fol})$.
  We now apply the $K$-theory machine  $\cK$ (see \cite[Def. 6.1]{2013arXiv1311.3188B} and Remark \ref{flwjfklefewfwfewfwfw234244}) and get a commutative diagram of presheaves of spectra  \begin{equation}\label{f4fwefwefewfwf}
\xymatrix{\cK(\Vect^{\flat,\nabla}) \ar[d]\ar[r]&\cK(\Vect^{\flat}) \ar[d]\\\cK(\Vect^{\nabla})\ar[r]\ar[d]&\cK(\Vect)\ar[d]\\   \widehat \ku^{\nabla}\ar[r]& \widehat \ku }
\end{equation}
  in $\PSh_{\Sp}(\Mf_{\C\textrm{-}\fol})$. 
   The upper square  in \eqref{f4fwefwefewfwf} is by definition the image  of \eqref{wsdqwdwqdqwdwqdsq} under $\cK$.
  The lower horizontal map is defined by applying the sheafification $L$ (see \eqref{dqdqwdwqdwqdwqdwqdqd}) to the middle horizontal arrow and the lower vertical arrows are the units of the sheafification.
  In particular, we use the notation \begin{equation}\label{wqdwqddq21321}
\widehat \ku^{\nabla}:=L(\cK(\Vect^{\nabla}))\ , \quad \widehat \ku:=L(\cK(\Vect))\ .
\end{equation}  
 
 \begin{rem}\label{flwjfklefewfwfewfwfw234244}{\rm 
For the sake of the reader let us indicate some details on the  $K$-theory machine $\cK$. It  is the composition \begin{eqnarray*}&&
 \CAlg(\Cat[W^{-1}])\to  \CAlg(\mathbf{Groupoids}[W^{-1}])\to \CommMon(\sSet[W^{-1}])\\&&\hspace{4cm}\to \CommGroup(\sSet[W^{-1}])\simeq \Sp_{\ge 0}\to \Sp \end{eqnarray*}

of the following functorial constructions:
\begin{enumerate}
\item We first take the underlying symmetric monoidal groupoid.
\item Then we apply the nerve in order to get a commutative monoid in  the category of spaces $\sSet[W^{-1}]$, i.e., an $E_{\infty}$-space.
\item Then we apply the group completion functor to obtain a commutative group in spaces, i.e., a grouplike $E_{\infty}$-space.
\item Finally we apply the functor which maps a commutative group in spaces to the corresponding    connective spectrum whose $\infty$-loop space is this group.
\end{enumerate}}
\end{rem}
 \begin{rem}{\rm 
 Note that the symmetric monoidal stacks $\Vect^{\nabla}$ and $\Vect$ are pulled back  
 from  stacks on the site $\Mf$ via the forgetful morphism \eqref{ewfwfwefewfewfwfwefwfe}. The same is true for the associated  sheaves of $K$-theory spectra $\widehat \ku^{\nabla}$ and $\widehat \ku$.  They represent  differential versions of connective $K$-theory $\ku$ and  are studied in detail in \cite[Sec. 6]{2013arXiv1311.3188B}
 }\end{rem}

 \subsection{Characteristic cocycles}\label{kfjwelfewfewfewfewfewfewfe}
 
 In order to construct the regulator we use the method introduced in  \cite{Bunke:2012fk} based on the notion of characteristic cocycles.   
 We consider the 
 category of chain complexes $\Ch$. 
 We have $$DD^{-}, DD^{\per}\in \Sh_{ \Ch }(\Mf_{\C\textrm{-}\fol})$$
 introduced in Definition \ref{jkdjlqwdqwdqwd}, where here we forget the algebra structure, and we use \cite[Lemma 7.12]{2013arXiv1311.3188B}  for the sheaf condition.
   
 Using the Chern character forms (Definitions \ref{ffwefwefewfewfwfw} and \ref{qldjqwldqwdqwdqwd}) and their naturality (equations \eqref{wqdqwdqwdwqdwqwqdqd} and \eqref{wqdqwdqwdwqdwqwqdqd1}) we   get   
characteristic cocycles (see  \cite[Def. 2.12]{Bunke:2012fk})
$$\ch^{-}:\pi_{0}(\Vect^{\flat,\nabla})\to Z^{0}(DD^{-}) \ , \quad 
 \ch:\pi_{0}(\Vect^{ \nabla})\to Z^{0}(DD^{\per}) \ .$$
 Here $\pi_{0}$ sends a symmetric monoidal category to its commutative monoid of isomorphism classes. 
 We will consider  commutative monoids as symmetric monoidal categories, see Remark \ref{dkwqnmqlkwdwqdwqdqd}.
 The following diagram in $\PSh_{\CAlg(\Cat[W^{-1}])}(\Mf_{\C\textrm{-}\fol})$ commutes:
\begin{equation}\label{kjwehbfjkfhekfehwfkewfef09890}
\xymatrix{\pi_{0}(\Vect^{\flat,\nabla})\ar[r]^{\ch^{-}}\ar[d]& Z^{0}(DD^{-})\ar[d]\\\pi_{0}(\Vect^{ \nabla})\ar[r]^{\ch}&Z^{0}(DD^{\per})}\ .  \end{equation}
We can now apply the algebraic $K$-theory machine described in Remark \ref{flwjfklefewfwfewfwfw234244}
and get the following commuting diagram  in $\PSh_{\Sp}(\Mf_{\C\textrm{-}\fol})$.
\begin{equation}\label{kjwehbfjkfhekfehwfkewfef0989012}
\xymatrix{\cK(\pi_{0}(\Vect^{\flat,\nabla})\ar[r]^{\cK(\ch^{-})}\ar[d]& \cK(Z^{0}(DD^{-}))\ar[d]\\\cK(\pi_{0}(\Vect^{ \nabla}))\ar[r]^{\cK(\ch)}&\cK(Z^{0}(DD^{\per}))}  \end{equation}

Let $$\mathbf{H}:\Ch[W^{-1}]\to \Sp$$ denote the Eilenberg-MacLane functor (see \cite[(22)]{Bunke:2012fk}). We will use the notation
 \begin{equation}\label{ewfewfewfewf444efewfewfewfwef}
\sigma^{\ge p}\DD^{-}:=\mathbf{H}(\sigma^{\ge p }DD^{-})\ , \quad \sigma^{\ge p}\DD^{\per}:=\mathbf{H}(\sigma^{\ge p}DD^{\per})
\end{equation} for $p\in \Z$, where 
$\sigma^{\ge p}$ (on the right-hand sides) is the so-called stupid truncation functor on chain complexes which sends a chain complex
$$\cdots \to C^{p-2} \to C^{p-1}\to C^{p}\to C^{p+1}\to C^{p+2}\to \cdots$$ to its  part 
$$\cdots \to 0\to 0\to C^{p}\to C^{p+1}\to C^{p+2}\to \cdots $$ of degree $\ge p$.
Note that 
$$\sigma^{\ge p}\DD^{-}\ , \sigma^{\ge p}\DD^{\per}\in  \Sh_{ \Sp}(\Mf_{\C\textrm{-}\fol})$$
by \cite[Lemma 7.12]{2013arXiv1311.3188B} and the fact that  the Eilenberg-MacLane functor $\mathbf{H}$ preserves limits, see also  \cite[Sec. 2.3]{Bunke:2014aa}.

The following construction is a case of  the general  construction of regulators  given in Definition \cite[Def. 2.14]{Bunke:2012fk}. 
If we consider an abelian group $A$ as a symmetric monoidal category, then we have a natural equivalence
 $\cK(A)\simeq \mathbf{H}(A)$, see \cite[Rem 2.13]{Bunke:2012fk}. Furthermore, for a chain complex $C$ we can view the inclusion $Z^{0}(C)\to C^{0}\to \sigma^{\ge 0}C$ as a natural morphism of chain complexes. Using these observations we can extend the diagram \eqref{kjwehbfjkfhekfehwfkewfef0989012} to the right and obtain the following  commuting diagram  in $\PSh_{\Sp}(\Mf_{\C\textrm{-}\fol})$.
 \begin{equation}\label{kjwehbfjkfhekfehwfkewfef098901555}
\xymatrix{\cK(\pi_{0}(\Vect^{\flat,\nabla}))\ar[r]^(0.53){\cK(\ch^{-})}\ar[d]& \cK(Z^{0}(DD^{-}))\ar[d]\ar[r]^{\simeq} &\mathbf{H}(Z^{0}(DD^{-}))\ar[d]\ar[r]& \sigma^{\ge 0}\DD^{-}\ar[d]\\\cK(\pi_{0}(\Vect^{ \nabla}))\ar[r]^{\cK(\ch)}&\cK(Z^{0}(DD^{\per}))\ar[r]^{\simeq}& \mathbf{H}(Z^{0}(DD^{\per}))\ar[r]& \sigma^{\ge 0}\DD^{\per} }  \end{equation}

We keep the outer square and extend it further to 
the commuting diagram 
\begin{equation}\label{kjwehbfjkfhekfehwfkewfef098901}
 \xymatrix{\cK(\Vect^{\flat,\nabla}) \ar[r]^(0.55){r(\ch^{-})}\ar[d]& \sigma^{\ge 0}\DD^{-}\ar[d]\\\ \cK(\Vect^{\nabla})\ar[d]^{u}\ar[r]^{r(\ch)} &\sigma^{\ge 0}\DD^{\per}\\\widehat \ku^{\nabla}\ar@{.>}[ur]_{r(\ch)}& }
\end{equation} in $\PSh_{\Sp}(\Mf_{\C\textrm{-}\fol})$.
In order to get the lower triangle we use that $\sigma^{\ge 0}\DD^{\per}$ is a sheaf and the universal property of the unit $u$ of the sheafification involved in the definition \eqref{wqdwqddq21321} of $\widehat \ku^{\nabla}$.

\

\subsection{The class $[V,\nabla^{I}]^{\alg}$}

 Let us fix a manifold $X$. We want to consider foliations  whose
 space of leaves is $X$. Trivial foliations of this type are obtained by taking the product of the typical leaf with $X$.
 In this way we actually obtain an inclusion of manifolds into foliations. More precisely we consider the
 functor   \begin{equation}\label{dwedewdewdwd}
j_{X}:\Mf \to \Mf_{\C\textrm{-}\fol}\ , \quad j_{X}(P):=(P,\cF_{\max})\times (X,\cF_{\min})\ .
\end{equation} A  manifold $Y$ also gives rise to endofunctors
  \begin{equation}\label{a1}
 i_{Y}:\Mf \to \Mf \ , \quad i_{Y} (P):= Y\times P\ ,
\end{equation}  and \begin{equation}\label{a2}
i_{Y}:\Mf_{ \C\textrm{-}\fol} \to \Mf_{\C\textrm{-}\fol}\ ,\quad  i_{Y} (P,\cF):=(Y,\cF_{\max})\times (P,\cF)\ .
\end{equation}    
    \bigskip
    
 The projection $Y\to *$ induces a morphism $\id\to i_{Y}^{*}$ on presheaves.

 \bigskip
 
 Let $I:=[0,1]$ denote the unit interval. Let $\bC$ be a presentable $\infty$-category. \begin{ddd}
An object $A\in \Sh_{\bC}(\Mf)$ (or $ \PSh_{\bC}(\Mf)$, $ \Sh_{\bC}(\Mf_{\C\textrm{-}\fol})$ or   $ \PSh_{\bC}(\Mf_{\C\textrm{-}\fol})$) is called homotopy invariant if the natural morphism
$A\to i_{I}^{*} A$ is an equivalence.  \end{ddd}
 We indicate the full subcategories of homotopy invariant (pre)sheaves by an upper index $h$.
 Note that in the foliated case we include the tangent bundle of the interval into the foliation direction.
 \bigskip
 
 \begin{ex}{\rm 
 By \cite[Prop. 2.6, 1.]{2013arXiv1311.3188B} (see also Lemma \ref{kfhkwjefewfewf98790234jnbkjf} below), for a homotopy invariant sheaf $\bE\in  \Sh^{h}_{\bC}(\Mf)$  we have a natural equivalence
 $$\underline{\bE(*)}\simeq \bE\ ,$$
 where    $\underline{\bE(*)}$ denotes the sheaf obtained from the constant presheaf with value $\bE(*)$ by sheafification. If $\bC=\Sp$, then for 
 $M\in \Mf$ and $k\in \Z$ we have a natural isomorphism of abelian groups
 \begin{equation}\label{fwefewfewfewfe324}
\pi_{k}(\bE(M))\cong \bE(*)^{-k}(M)\ .
\end{equation} 
Observe that a similar statement is not true for homotopy invariant sheaves on $\Mf_{\C\textrm{-}\fol}$.}
\end{ex}

 \begin{lem}\label{leelfjwelfwefewf}
The sheaf $\Vect^{\flat}$ is homotopy invariant. 
\end{lem}
\proof
The reason is that the foliation of $i_{I}(M,\cF)=(I\times M,T_{\C}I\boxplus \cF)$ contains the $I$-direction. For $(V,\nabla^{I})\in \Vect^{\flat}(i_{X}(M,\cF))$ 
we can use the flat connection $\nabla^{I}$ in order to define a parallel transport in the $I$-direction. 

Hence a vector bundle $(V,\nabla^{I})$ with a flat partial connection or a morphism between two such objects over $I\times M$ is uniquely determined by the restriction to $\{0\}\times M$. \hB

We now use the fact that on the site $\Mf$ sheafification preserves homotopy invariance.
\begin{lem}\label{kfhkwjefewfewf98790234jnbkjf}
If $\bbF\in \PSh^{h}_{\bC}(M)$, then $L(\bbF)\simeq \underline{\bbF(*)} $. In particular $L(\bbF)\in \Sh^{h}_{\bC}(\Mf)$.
\end{lem}
\proof
This is \cite[Prop. 2.6.2.]{2013arXiv1311.3188B}
\hB

We define \begin{equation}\label{vrekvejknrvkjenvkjvrevrrev}
\bK_{X}:=L(j_{X}^{*} \cK(\Vect^{\flat}))\in \Sh^{h}_{\Sp}(\Mf)\ .
\end{equation} 
 Note that $j_{X}^{*}$ preserves homotopy invariance and the sheaf condition. By Lemmas \ref{leelfjwelfwefewf}
  and \ref{kfhkwjefewfewf98790234jnbkjf} we see that $\bK_{X}$ is indeed a homotopy invariant sheaf.  

\bigskip

 We have a chain of equivalences of symmetric monoidal categories $$j_{X}^{*}\Vect^{\flat}(*)\simeq \Vect^{\flat}(X,\cF_{\min})\simeq \Vect(X)\simeq \Proj(C^{\infty}(X))\ ,$$
 where the first three are obtained by specializing definitions, and the last is Swan's theorem.
This implies 
 \begin{equation}\label{fefwefwefewf234234sdfsdf}
\bK_{X}(*)= \cK(\Proj(C^{\infty}(X)))\stackrel{\mathrm{def}}{=}\bK(C^{\infty}(X))\ ,
\end{equation}
where the last equality is our definition of the connective algebraic $K$-theory spectrum of the ring $C^{\infty}(X)$.
 
We can now give the technical definition of the class $[V,\nabla^{I}]^{\alg}\in \bK(C^{\infty}(X))^{0}(P)$ for a pair $$(V,\nabla^{I})\in \Vect^{\flat}(P\times X,T_{\C}P\oplus 0)\ .$$ Indeed, we have
$(V,\nabla^{I})\in j_{X}^{*}\Vect^{\flat}(P)$. This object naturally represents a  point\footnote{Note that an object in a symmetric monoidal category $\bC$ naturally represents a point in the nerve $\Nerve(\bC)$ of $\bC$ in therefore a point (up to contractible choice, i.e., a component) in its group completion. The latter    is, by definition, $\Omega^{\infty}\cK(\bC)$, see  \ref{flwjfklefewfwfewfwfw234244}.
} in $\Omega^{\infty} \bK_{X}(P)$. 
\begin{ddd}\label{xgrergegg} We define  \begin{equation}\label{fwefwefwefffefefewfewfewff324234}
[V,\nabla^{I}]^{\alg}\in \pi_{0}(\bK_{X}(P))\stackrel{\eqref{fwefewfewfewfe324}}{\cong} \bK_{X}(*)^{0}(P) \stackrel{\eqref{fefwefwefewf234234sdfsdf}}{\cong} \bK(C^{\infty}(X))^{0}(P)\ .
\end{equation} 
to be the connected component represented by the point $[V,\nabla^{I}]$.\end{ddd}
 
\subsection{Differential $K$-theory and the regulator map}

We assume that $\bC$ is a stable presentable $\infty$-category like spectra $\Sp$ or chain complexes $\Ch[W^{-1}]$.
We have   an adjunction
$$\cH:\Sh_{\bC}(\Mf  )\leftrightarrows \Sh^{h}_{\bC}(\Mf ):\mathrm{inclusion}\ ,$$
where $\cH$ is called the homotopification. By  \cite[Prop. 7.6.(2)]{2013arXiv1311.3188B} it is given by a composition
$\cH\simeq L\circ \cH^{\mathrm{pre}}$, where $\cH^{\mathrm{pre}}:\Sh_{\bC}(\Mf  )\to \PSh^{h}_{\bC}(\Mf)$ is given by 
\begin{equation}\label{fewfwfewfewfewfewfewfewfewfewf}
\cH^{\mathrm{pre}} \simeq \colim_{\Delta^{\op}}  i_{\Delta^{\bullet}}^{*} 
\end{equation}
using the notation \eqref{a1}.
Similarly, for the site $\Mf_{\C\textrm{-}\fol}$ we have an adjunction \begin{equation}\label{ggg889899898893443}
\cH^{\flat}:\Sh_{\bC}(\Mf_{\C\textrm{-}\fol}  )\leftrightarrows \Sh^{h}_{\bC}(\Mf_{\C\textrm{-}\fol} ):\mathrm{inclusion}\  ,
\end{equation} where $\cH^{\flat}=L\circ \cH^{\flat,\mathrm{pre}}$ with $\cH^{\flat, \mathrm{pre}}$ given again by \eqref{fewfwfewfewfewfewfewfewfewfewf}, but now using \eqref{a2}.
 For a manifold $X$ the functor $j_{X}^{*}$ (see \eqref{dwedewdewdwd}) preserves homotopy invariant sheaves. Moreover, if $X$ is compact, then we have \begin{equation}\label{refwefwefewfew45355321343241325465}
j_{X}^{*}\circ \cH^{\flat}\simeq \cH\circ j_{X}^{*}
\end{equation}
(compare \cite[Lemma 2.4 (4)]{Bunke:2014aa} for a proof of  a similar statement).

\begin{lem} \label{iewfwefewfewfewf}The sheaves $\DD^{\per}$ and $\DD^{-}$ are homotopy invariant. Moreover, for every $p\in \Z$ the inclusions
$$\sigma^{\ge p} \DD^{-}\to \DD^{-}\ , \quad \sigma^{\ge p} \DD^{\per}\to \DD^{\per}$$ are equivalent to the units of the homotopification.
\end{lem} \proof We start with the case of the map $\sigma^{\ge p} \DD^{\per}\to \DD^{\per}$ between sheaves on $\Mf$. 
Recall the definition \eqref{ewfewfewfewf444efewfewfewfwef}. We let $$\iota:\Ch\to \Ch[W^{-1}]$$ be the canonical localization map. 
We have $$\sigma^{\ge p}DD^{\per}\cong \prod_{q\in \Z} (\sigma^{\ge p+2q}\Omega)[2q]\ .$$
We discuss the factors separately.
By \cite[Lemma 7.15]{2013arXiv1311.3188B} the map
$$\iota (\sigma^{\ge p+2q}\Omega)[2q] \to \iota( \Omega)[2q]$$ is the unit of the homotopification.
This implies the assertion for $ \DD^{\per}$ after applying the Eilenberg-MacLane functor $\mathbf{H}$.

\bigskip

We now discuss $DD^{-}$. We first observe that
$\iota (DD^{-})$ is a homotopy invariant sheaf  on the site $\Mf_{\C\textrm{-}\fol}$ with values in $\Ch[W^{-1}]$.
We again consider one factor   of $$ DD^{-}\cong \prod_{q\in \Z}  F^{q}\Omega [2q]$$ at a time. For a foliated manifold $(M,\cF)$ the integration $\int_{I\times M/M}$ preserves the filtration and induces a map
$$\int_{I\times M/M} \colon F^{q}\Omega(I\times M,T_{\C}I\boxplus \cF) \to F^{q}\Omega(M, \cF)[ -1]$$ such that
$$d\int_{I\times M/M} x=x_{|\{1\}\times M}-x_{|\{0\}\times M}\ .$$
This implies that $\iota (F^{q}\Omega )$ is homotopy invariant.

\begin{rem}{\rm 
The point here is that we define homotopy invariance along the leaf direction.
If we would include transverse directions, then the integral would not preserve the filtration. In this case we only have 
$$\int_{I\times M/M} F^{p}\Omega(I\times M,  \{0\}\boxplus\cF) \to F^{p-1}\Omega(M, \cF)[ -1]$$
and  the integration would not be defined on $DD^{-}$. 
 }
\end{rem}

Once we know that $\iota (F^{q}\Omega) \in \Sh_{\Ch[W^{-1}]}(\Mf_{\C\textrm{-}\fol})$ is homotopy invariant, we show 
that $$\iota( \sigma^{\ge p} F^{q}\Omega) \to  \iota  (F^{q}\Omega) $$
is the unit of the homotopification exactly as in \cite[Lemma 7.15]{2013arXiv1311.3188B}.
Note that by (the analogue of) \cite[Lemma 7.13]{2013arXiv1311.3188B}  $\cH^{\flat}(\iota (F^{q}\Omega)^{\ell})=0$ for every $\ell\in \Z$.
This implies as in the proof of \cite[Lemma 7.15]{2013arXiv1311.3188B} that
$\cH^{\flat}(\iota (\sigma^{<p} F^{q}\Omega) )=0$. The claim now follows from an application of $\cH^{\flat}\circ \iota$ to the exact sequence of $\Ch$-valued sheaves
$$0\to\sigma^{<p} F^{q}\Omega \to   F^{q}\Omega  \to \sigma^{\ge p} F^{q}\Omega  \to 0\ . $$

\hB

 Recall the definition \eqref{vrekvejknrvkjenvkjvrevrrev} of  the sheaf of spectra $\bK_{X}$.
We define $$\bK^{\nabla}_{X}:=L(j_{X}^{*} \cK(\Vect^{\flat,\nabla}))\ .$$

\begin{lem}\label{kljldqwdqwdwqdd}
 The morphisms $$ \bK_{X}^{\nabla}\to \bK_{X}\ , \quad \widehat \ku \to \underline{\ku}\ , \quad \widehat \ku^{\nabla}\to \underline{\ku}$$
are equivalent to  the units of the homotopification.
\end{lem}
\proof
The second and the third cases are consequences of \cite[Lemma 6.3]{2013arXiv1311.3188B}
and \cite[Lemma 6.5]{2013arXiv1311.3188B}. It remains to discuss the first case.  
We know that $\bK_{X}$ is homotopy invariant.  Then the assertion   follows from the analogue of
\cite[Lemma 6.4]{2013arXiv1311.3188B} for $\Vect^{\flat,\nabla}\to\Vect^{\flat}$. \hB

 From \eqref{kjwehbfjkfhekfehwfkewfef098901} and the fact that the two objects on the right and the lower left corner are sheaves 
 we get the diagram 
\begin{equation}\label{kjwehbfjkfhekfehwfkewfeeedef098901}
 \xymatrix{\bK^{\nabla}_{X}  \ar[rr]^{j_{X}^{*}r(\ch^{-})}\ar[d]&& j_{X}^{*}\sigma^{\ge 0}\DD^{-}\ar[d]\\ j_{X}^{*}\widehat \ku^{ \nabla}\ar[rr]^{j_{X}^{*}r(\ch)}&&j_{X}^{*}\sigma^{\ge 0}\DD^{\per}}\ .
\end{equation}
We now assume that $X$ is compact. Then by \eqref{refwefwefewfew45355321343241325465}
 homotopification commutes with $j_{X}^{*}$. 
Applying homotopification to this square and using Lemmas \ref{iewfwefewfewfewf}, \ref{kljldqwdqwdwqdd} we get the square  \begin{equation}\label{r223r23r32r32r324}
 \xymatrix{\bK_{X} \ar[r]^{\omega^{-}_{X}}\ar[d]& j_{X}^{*} \DD^{-}\ar[d]\\ j_{X}^{*}\underline{\ku} \ar[r]^{\omega_{X}}& j_{X}^{*} \DD^{\per}}\ . \end{equation}

We consider the following three versions of Hopkins-Singer type (see \cite{MR2192936} for the original definition and \cite{2013arXiv1311.3188B} for more information) differential algebraic and differential $K$-theories for $p\in \Z$
$$\xymatrix{\hat \bK_{X }^{p} \ar[d]^{I}\ar[r]&j_{X}^{*}\sigma^{\ge p}\DD^{-}\ar[d]\\ \bK_{X } \ar[r]^{\omega^{-}_{X}}&j^{*}_{X}\DD^{-}\ ,} \quad \xymatrix{\widehat \ku_{X}^{\flat,p}\ar[d]\ar[r]^{R}&j_{X}^{*}\sigma^{\ge p}\DD^{-}\ar[d]\\ j_{X}^{*}\underline{\ku}\ar[r]^{\omega_{X}}&j_{X}^{*}\DD^{\per}\ ,}
 \quad 
 \xymatrix{\widehat \ku^{p}\ar[d]\ar[r]&\sigma^{\ge p}\DD^{\per}\ar[d]\\\underline{\ku}\ar[r]&\DD^{\per}}$$
 defined by the respective pull-back square in $\Sh_{\Sp}(\Mf)$. We define the corresponding differential cohomology groups by
 $$\hat K_{X}^{p}(P):=\pi_{-p}(\hat \bK_{X }^{p} (P))\ , \quad  \widehat{ku}_{X}^{\flat,p}(P):=\pi_{-p}(\widehat \ku_{X }^{\flat, p} (P))\ , \quad 
  \widehat{ku}^{p}(P):=\pi_{-p}(\widehat \ku^{p}(P))\ .$$

The square \eqref{r223r23r32r32r324} together with the obvious commutative square  \begin{equation}\label{hjkdekdhewkdjhewkdewd}
\xymatrix{\sigma^{\ge p}\DD^{-} \ar[r]\ar[d]&  \DD^{-}\ar[d]\\ \sigma^{\ge p}\DD^{\per} \ar[r]&  \DD^{\per}}\ .\end{equation}
  induces a chain of morphisms
$$\hat \bK_{X}^{p}\to \widehat \ku_{X}^{\flat,p}\to  j^{*}_{X}\widehat \ku^{p}\ .$$ 
Using \eqref{kjwehbfjkfhekfehwfkewfeeedef098901} we finally get the square  
 $$\xymatrix{  \bK_{X}^{\nabla} \ar[r]^{\mathrm{cycl}}\ar[d]& \hat \bK_{X}^{0}\ar[d]\\ j_{X}^{*}\ku^{\nabla} \ar[r]^{\mathrm{cycl}}&j_{X}^{*}\widehat \ku^{0}}$$ where the horizontal maps are the differential cycle maps.

\bigskip

The following exact sequences are part of the general features of a Hopkins-Singer differential cohomology, see e.g. \cite[Proposition 2.3.1.]{Bunke:2013ab} or \cite[Remark 4.9]{2013arXiv1311.3188B}. The sequence 
\begin{equation}\label{}   \dots\to    DD^{-}(P\times X,T_{\C}P\boxplus\{0\})^{\ell-1} /\im(d) \stackrel{a}{\to} \hat K_{X}^{\ell}(P)\stackrel{I}{\to}  \bK_{X}(*)^{\ell}(P)\to 0
\end{equation}
describes the set of possible differential lifts of topological classes. The second sequence
\begin{equation}\label{} 0\to \widehat{ku}^{\flat,\ell}_{X,\flat}(P)\to \widehat{ku}^{\flat,\ell}_{X}(P)\stackrel{R}{\to} Z^{\ell}(DD^{-}(X\times P,T_{\C}P\boxplus\{0\}))\to\dots
\end{equation}
reflects the definition of the flat subgroup.

\bigskip

 We consider the case $P=*$, $\ell:=-p$ and assume that $\dim(X)<p$. Then it is straightforward to check that  $DD^{-}(X,\cF_{\min})^{-p-1}=0$ and $DD^{-}(X,\cF_{\min})^{-p}=0$.
 

This implies the isomorphisms $$I:\hat K_{X}^{-p}(*)\stackrel{\cong}{\to} \bK_{X}(*)^{-p}\ , \quad  \widehat{ku}_{X}^{\flat,-p}(*)\stackrel{\cong}{\to} \widehat{ku}_{X,\flat}^{\flat,-p}(*)\ .$$  
\begin{ddd}\label{klfwefewfewfwf} For $p\in \nat$ such that $\dim(X)<p$ we define the regulator map
  $\reg_{X}$ as the composition
$$\hspace{-1cm}K_{p}(C^{\infty}(X))\cong \bK_{X}(*)^{-p}\stackrel{\cong}{\leftarrow} \hat K_{X}^{-p}(*) \to \widehat{ku}_{X}^{\flat,-p} (*)\stackrel{\cong}{\rightarrow} \widehat{ku}^{\flat,-p}_{X,\flat} (*)
\to \widehat{ku}^{-p}_{\flat}(X,\cF_{\min})\stackrel{!}{\cong} \ku\C/\Z^{-p-1}(X)$$
 \end{ddd} 
For the natural isomorphism marked by $!$ in the formula above we refer to \cite[Proposition 2.3.2.]{Bunke:2013ab} 
or \cite[Remark 4.9]{2013arXiv1311.3188B}.
In  Remark \ref{2ddhi3dhio32doi2doo2oidjud2} we will explain how this regulator can be obtained by specializing a more basic regulator.

\begin{rem}{\rm 
In \cite[Thm 1.1]{Bunke:2014aa} we defined a similar regulator map 
$$\sigma_{p}:K_{p}(C^{\infty}(X))\to  \ku\C/\Z^{-p-1}(X)$$
using different methods. While here, in order to define the Chern character, we use characteristic forms associated to connections, in  \cite{Bunke:2014aa} we use the  Goodwillie-Jones Chern character. 
The two Chern characters are equivalent as primary invariants \cite[Lemma 2.27]{Bunke:2014aa}. In order to compare
the two regulator maps $\sigma_{p}$ and $\reg_{X}$ we would need to compare the two Chern characters on the space level. So at the moment it remains an open question whether $\sigma_{p}=\reg_{X}$.
}\end{rem}



\subsection{Integration and proof of Theorem \ref{flkfefwefwefewfef}}\label{fwklfwfewfewfewf}

We now assume that $P$ is closed and  has a stable framing $s$. Then $f:P\to *$ has a natural differential orientation $\hat o_{s}$  (see  \cite[Example 4.230]{2012arXiv1208.3961B})  and we have an associated Umkehr map in every Hopkins-Singer differential cohomology theory.
We further assume that $X$ is closed, spin and equipped with a Riemannian metric. This induces a differential $\ku$-orientation $\hat o$ of the projection $\pi:X\to *$, see Subsection \ref{wlekfjwelfjewlfwef123}.

Let  $p:=\dim(P)$ and $d:=\dim(X)$. We have the commutative diagram 
\begin{equation}\label{dewewdedewdewd2342343}
\xymatrix{ &&\widehat{KU}^{0}(P\times X)\ar[r]^{(\pi\circ f)^{\hat o\circ \hat o_{s}}_{!}}&\widehat{KU}^{-p-d}(*)\ar[r]^{\cong}&\C/\Z\ar@{=}[d]\\
\pi_{0}(\bK_{X}^{\nabla})\ar@/^3cm/[rrrru]^{(V,\nabla^{V})\mapsto \rho(M,\cF,\nabla^{I},s)} \ar[urr]^{(V,\nabla)\mapsto [V,\nabla]}\ar[d]_{(V,\nabla^{V})\mapsto [V,\nabla^{I}]^{\alg}}\ar[r]^{\mathrm{cycl}} &\hat K_{X}^{0}(P)\ar[r]\ar[dl]\ar[d]^{f_{!}^{\hat o_{s}}}\ar@{}[dr]^{\textcircled{2}}&\widehat{ku}^{0}(P\times X) \ar@{}[ur]^{\textcircled{2}}\ar[d]^{\hat f^{\hat o_{s}\times X}_{!}}\ar@{}[dr]^{\textcircled{1}}\ar[u]\ar[r]^{(\pi\circ f)^{\hat o\circ \hat o_{s}}_{!}}&\widehat{ku}^{-p-d}(*) \ar[u]\ar@{=}[d]\ar[r]^{\cong}&\C/\Z\ar@{=}[d]\\
   \bK_{X}(*)^{0}(P)\ar@{}[r]^{\textcircled{3}} \ar[d]^{f_{!}^{o_{s}}}&\ar[dl]^{\cong}\hat K^{-p}_{X}( *)  \ar@{}[dr]^{\textcircled{4}}\ar[r] &\widehat{ku}^{-p }(X)     \ar[r]^(0.5){\pi_{!}^{\hat o}}  &\widehat{ku}^{-d-p}(*) \ar[r]^{\cong}&\C/\Z 
   \\ \bK_{X}(*)^{-p} \ar[rr]^{\reg_{X}}
 & &\ku\C/\Z^{-p-1}(X)\ar[r]^{\pi_{!}^{o}} \ar[u]^{!!}\ar@{}[ur]^{\textcircled{5}}&\ku\C/\Z^{-d-p-1}(*)\ar[u]_{\cong } & }
\end{equation}
 The   square $\textcircled{1}$  commutes by the $\ku$-analogue of  \eqref{e23e23e3e23e32e32e2}.
For the squares $\textcircled{2}$  we use that integration commutes with transformations between Hopkins-Singer differential cohomology theories provided the orientations are related correspondingly.
For the  square $\textcircled{3}$ we use the right-most square of  the $\ku$-analogue of \eqref{r23r23r23r23r23r235435346546}.
The square $\textcircled{4}$ commutes by the Definition \ref{klfwefewfewfwf} of the regulator. Here we also use the square 

$$\xymatrix{ \ku\C/\Z^{-p-1}(X)   \ar[r]^{!!} & \widehat{ku}^{-p}(X)  \\  \widehat{ku}^{-p}_{\flat}(X,\cF_{\min}) \ar[r] \ar[u]_{\cong}& \widehat{ku}^{-p}_{\flat} (X) \ar[u]_{\subseteq}}$$ explaining the arrows marked by $!!$   above and in 
\eqref{dewewdedewdewd2342343}.
For $\textcircled{5}$
 we use that the identification of the flat subgroup in a Hopkins-Singer differential cohomology with  the $\C/\Z$-version of the underlying cohomology theory is compatible with integration, i.e., the left-most square in  the $\ku$-analogue of \eqref{r23r23r23r23r23r235435346546}.

\bigskip

The upper composition in \eqref{dewewdedewdewd2342343} maps $(V,\nabla^{V})$, essentially by definition,  to $\rho(M,\cF,\nabla^{I},s)$ as indicated. The down-right composition sends $(V,\nabla^{V})$ to
$$\pi_{!}^{o}(\reg_{X}(f_{!}^{o_{s}}([V,\nabla^{I}]^{\alg})))\ .$$
Thus Theorem \ref{flkfefwefwefewfef} follows from the commutativity of \eqref{dewewdedewdewd2342343}. \hB
 
\section{Algebraic $K$-theory of foliations}\label{dkqwldqwdwqdwqdwqdwqd}

In this section we define the algebraic $K$-theory sheaf $\bK$ on $\Mf_{\C\textrm{-}\fol}$. Its homotopy groups $$K^{*}(M,\cF):=\pi_{-*}(\bK(M,\cF)) $$ can be considered as the algebraic $K$-theory groups of the foliation $(M,\cF)$. We further introduce the Hodge-filtered connective $K$-theory sheaf $  \ku^{\flat}$ and define a regulator
$$\reg:\bK\to  \ku^{\flat}\ .$$ For $p>\codim(\cF)$ it induces a map
$$\widetilde \reg:K^{-p}(M,\cF)\to \ku\C/\Z^{-p-1}(M)$$
 which generalizes the regulator introduced in Definition \ref{klfwefewfewfwf}.

\bigskip

\begin{rem}{\rm 
This section has a considerable overlap with the work of Karoubi \cite{karoubi43}, \cite{karoubi45}.
We add this section to the present paper since it fits well with the set-up developed here and puts the regulator in its natural framework. 
We will study this regulator and examples elsewhere.  
}
\end{rem}

We will use the notation introduced in Subsection \ref{klfwjlefewfewf}.  In particular, $\Vect^{\flat}$ and $\Vect^{\flat,\nabla}$ denote the symmetric monoidal stacks of pairs $(V,\nabla^{I})$ and $(V,\nabla)$ of complex vector bundles and flat partial connections, or complex vector bundles and connections whose restriction to the foliation is flat, respectively.
The symbols $L$ and $\cH^{\flat}$ denote the sheafification and the homotopification operations.

 \begin{ddd} We define
  sheaves of spectra
$$ \bK :=  \cH^{\flat}(L(\cK(\Vect^{\flat}))) \in \Sh^{h}_{\Sp}(\Mf_{\C\textrm{-}\fol})\ , \quad  \bK^{\nabla} :=  L(\cK(\Vect^{\flat,\nabla})) \in \Sh_{\Sp}(\Mf_{\C\textrm{-}\fol})\ .$$ For $p\in \Z$ we 
define the algebraic $K$-theory of a foliated manifold $(M,\cF)$ by
$$ K^{p}(M,\cF):=\pi_{-p}( \bK (M,\cF))\ .$$
 \end{ddd}
\begin{rem}\rm {Note that $\cK(\Vect^{\flat})$ is homotopy invariant. 
We expect that the sheafification preserves homotopy invariance so that the homotopification is not really necessary in this definition.
}\end{rem}

In order to motivate this definition let us discuss special cases.

\begin{ex}{\rm 
Recall the functor $j:=j_{*}:\Mf\to \Mf_{\C\textrm{-}\fol}$ given by $$j(M):=(M,\cF_{\max})\ ,$$ see \eqref{dwedewdewdwd}. Let  $\bK(\C)$ denote  the  connective algebraic $K$-theory spectrum of the field $\C$.
\begin{lem} \label{fwekflewfwfewfw} We have an equivalence $j^{*}\bK\simeq \underline{\bK(\C)}$ \end{lem} \proof
Since 
 $j (I\times M)\cong (I, T_{\C}I)\times j(M)$ we 
 conclude that $j^{*}$ preserves homotopy invariant sheaves. Since $\bK$ is homotopy invariant the sheaf
 $j^{*}\bK$ is homotopy invariant. Therefore (see \cite[Prop. 2.6, 1.]{2013arXiv1311.3188B}) we have an equivalence
 $$j^{*}\bK\simeq \underline{(j^{*}\bK)(*)}\ .$$
 If $\bE$ is a presheaf of spectra on $\Mf_{\C\textrm{-}\fol}$ and $L$ is the sheafification \eqref{dqdqwdwqdwqdwqdwqdqd}, then we have a natural equivalence of spectra $L(\bE)(*)\simeq \bE(*)$.
 Consequently, 
 $$(j^{*}\bK)(*)\simeq \cK(\Vect^{\flat}(*,\cF_{\max}))\ .$$  The category $\Vect^{\flat}(*,\cF_{\max}) $ is  the category of finite-dimensional complex vector spaces. Consequently    we have an equivalence of spectra $$\cK(\Vect^{\flat}(*,\cF_{\max}))\simeq \bK(\C)\ .$$ 
 The combination of these equivalences gives the assertion of the lemma.
 \hB
 
 \bigskip
 
As a consequence of  Lemma \ref{fwekflewfwfewfw} we have for a manifold $M$  \begin{equation}\label{cejkwdjejdhewdhewkdei832e32e32e32e2e}
K^{*}(M,\cF_{\max})\cong \bK(\C)^{*}(M)\ .
\end{equation}

%
%
%
}
\end{ex}

\begin{ex}\label{kjffjewfewjkfhewkjfewfkewfhewf87z}{\rm 
We have a natural functor
 $\kappa:\Mf\to \Mf_{\C\textrm{-}\fol}$ which is given by $M\mapsto (M,\cF_{\min})$.
 On the site $\Mf$ we have the differential cohomology theory $\widehat \ku$, see \eqref{wqdwqddq21321} and \cite{2013arXiv1311.3188B}.
 We have an equivalence of sheaves of spectra on $\Mf_{\C\textrm{-}\fol}$ $$\kappa^{*}\bK\simeq \widehat \ku\ .$$ 
Consequently, $$K^{*}(M,\cF_{\min})\cong \widehat{ku}^{*}(M)\ .$$


 }
\end{ex}

\begin{ex}{\rm 
For a fixed manifold $P$ there is a natural map \begin{equation}\label{dkejwdlewdewded}
\bK_{X}(P)\to \bK(P\times X,T_{\C}P\boxplus\{0\})\end{equation}
which is natural in $X$. It
  is essentially the sheafification morphism in the direction of $X$.
 The spectrum valued functor $$P\mapsto \bK(P\times X, T_{\C}P\boxplus 0)$$ is a homotopy invariant sheaf on $\Mf$. 
 Since
 $$\bK(\{*\}\times X, T_{\C}\{0\}\boxplus 0)\simeq \widehat \ku(X)$$ by  \cite[Prop.2.6, 1.]{2013arXiv1311.3188B}) it
  is therefore equivalent to $\underline{\widehat \ku(X)}$. We thus  get a map 
$$  K_{X}^{*}(P)\to K^{*}(P\times X, T_{\C}P\boxplus 0)\cong \widehat{\ku}(X)^{*}(P)\ .$$
}
\end{ex}

\begin{ex}\label{kdjqwkldjwqldjwqldwqd}{\rm 
Assume that $X$ is a smooth complex algebraic variety and let $X^{\an}$ be its 
  associated complex manifold with the foliation $\cF:=T^{0,1}X^{\an}$.  Then we can consider the algebraic $K$-theory of $\bK^{\alg}(X)$. It is defined like $\bK(M,\cF)$ as the sheafification of the presheaf $X\supset U\mapsto \cK(\Vect^{\alg}(U))$, where $\Vect^{\alg}(U)$ is the symmetric monoidal category of algebraic vector bundles on the Zariski open subset $U$. Since the analytic topology of $M$ refines the Zariski topology of $X$ the
transformations  $\Vect^{\alg}(U)\to \Vect(U^{\an})$ induce a map
$$\bK^{\alg}(X)\to \bK(X^{\an},T^{0,1}X)\ .$$
This example justifies to call $\bK(M,\cF)$ the algebraic $K$-theory spectrum of the foliated manifold $(M,\cF)$.
    }
\end{ex}

\begin{ex}{\rm 
If $(V,\nabla^{I})$ is a complex vector bundle with flat partial connection on a foliated manifold $(M,\cF)$, then we get a class
$$[V,\nabla^{I}]^{\alg}\in K^{0}(M,\cF)\ .$$
Similarly, if $\nabla$ is a connection which extends $\nabla^{I}$, then we get a class
$$[V,\nabla]^{\alg}\in  \pi_{0}(\bK^{\nabla}(M,\cF))\ .$$}
\end{ex}

\bigskip

 From \eqref{kjwehbfjkfhekfehwfkewfef098901} and the fact that the objects on the right and the lower left corner are sheaves 
 we get the diagram 
\begin{equation}\label{kjwehbfjkfhekfehwfkewfeeedef098901rrr}
 \xymatrix{\bK^{\nabla}  \ar[r]^{ r(\ch^{-})}\ar[d]&  \sigma^{\ge 0}\DD^{-}\ar[d]\\  \widehat \ku^{ \nabla}\ar[r]^{ r(\ch)}& \sigma^{\ge 0}\DD^{\per}}\ .
\end{equation}
 Applying homotopification to this square and using the Lemmas \ref{iewfwefewfewfewf} and  \ref{kljldqwdqwdwqdd} 
 we get   the square  \begin{equation}\label{r223r23r32r32r324rrr}
 \xymatrix{\bK  \ar[r]^{\omega^{-}}\ar[d]&  \DD^{-}\ar[d]\\  \underline{\ku} \ar[r]^{\omega}&  \DD^{\per}}\ . \end{equation}

 \begin{ddd}
We define the Hodge-filtered connective complex $\ku$-theory sheaf $\ku^{\flat}$ on $\Mf_{\C\textrm{-}\fol}$ by the pull-back square \begin{equation}\label{ffewfewfwefewfewf234}
 \xymatrix{  \ku^{\flat}\ar[d]\ar[r]&\DD^{-}\ar[d]\\
\underline{\ku}\ar[r]^{\omega}&\DD^{\per}}\ .
\end{equation}
We let $$ ku^{\flat,p}(M,\cF):=\pi_{-p}(  \ku^{\flat}(M,\cF))$$ be the corresponding Hodge-filtered  $\ku$-theory groups of $(M,\cF)$.
\end{ddd}

\begin{rem}{\rm 
In \cite{karoubi43}, \cite{karoubi45}   Karoubi introduced, starting from a filtration of the de Rham complex, the  multiplicative K-theory $\mathbf{MK}$.  Applied to the
filtration (Definition \ref{ilfjewlfwfewfewfewfwfw}) coming from a foliation the multiplicative
K-theory groups $\mathbf{MK}^{*}(M,\cF)$ are the Hodge-filtered $\KU$-theory groups of $(M,\cF)$. In other words, $ku^{\flat,*}$ is the connective $K$-theory analogue of Karoubi's
multiplicative $K$-theory. If one applies the functor $\Omega^{\infty}$ to \eqref{ffewfewfwefewfewf234}, then one obtains a
pull-back square of sheaves of spaces which is the analogue of the square just before the statement of Theorem 7.3 in \cite{karoubi45}. The fact that $\ku^{\flat}$ is a sheaf of spectra implies a Mayer-Vietoris type sequence for an open decomposition of a foliated manifold. This is Karoubi's theorem  \cite[Thm. 7.7]{karoubi45}.

For a justification to use the term {\em Hodge-filtered...}  instead of {\em multiplicative...}
see Remark \ref{jhdjkhkjqwdhqwkjdhwqkjdqwdwqd}.

}
\end{rem}

\begin{rem}\label{jhdjkhkjqwdhqwkjdhwqkjdqwdwqd}{\rm The Hodge-filtered connective complex $\ku$-theory $\ku^{\flat}$ is the $\ku$-theory analogue of the integral Deligne cohomology which would be the Hodge filtered version of $H\Z$. While 
integral Deligne cohomology is the natural target for cycle maps from Chow groups of algebraic cycles,  $\ku^{\flat}$ is the natural target of the regulator from algebraic $K$-theory.  In  \cite{MR3335251} the authors defined for  every spectrum over $\mathbf{H} \Z$ a  Hodge-filtered version. In an analogous manner, replacing integral Deligne cohomology by $\ku^{\flat}$ one could construct Hodge filtered cohomology theories for spectra over $\ku$. Observe that
$\ku^{\flat}$ is the Hodge-filtered version associated to the identity $\ku\to \ku$. This fact motivates the name.
}
\end{rem}
In view of Lemma \ref{iewfwefewfewfewf} the sheaf 
$ \ku^{\flat}$ is homotopy invariant (compare \cite[Thm 4.8]{karoubi43}). This fact is reflected in our notation  by not   using a $\hat{(\dots)}$-decoration.

\bigskip

\begin{ddd}\label{ojlfewfefwefw234234}
We define the regulator $\reg:\bK\to \ku^{\flat}$ to be the morphism induced by the square \eqref{r223r23r32r32r324rrr} and the universal property of the pull-back square \eqref{ffewfewfwefewfewf234}.
\end{ddd}

\begin{rem}
{\rm Such a regulator has first been defined in \cite[Sec.4]{karoubi45}. Karoubi's regulator provides a factorization
$$\bK\to \mathbf{MK}\to \bKU$$
of the map from algebraic to topological $K$-theory.
Our analogue is
$$\bK\to \ku^{\flat}\to \ku\ .$$

}
\end{rem}

\bigskip

\begin{rem}{\rm
The map
$$\reg:\bK\to \ku^{\flat}$$ could be considered as a foliated and integral analogue of   Beilinson's regulator.
In order to see this we show that the classical Beilinson regulator can be factored over the regulator $\reg$ defined above.

We first interpret real Deligne cohomology as Hodge-filtered $\ku\R$-theory. Here as usual,   we write $\ku\R:=\ku\wedge \mathbf{M}\R$ for the product  of $\ku$ with the Moore spectrum of $\R$.
The Chern character induces an equivalence of spectra $$\ku\R\simeq \prod_{p\ge 0}\mathbf{H}\R[2p]\ .$$
The de Rham equivalence $\underline{\mathbf{H}\R}\simeq \mathbf{H}\Omega_{\R}$ provides  the second  equivalence in the composition
$$\underline{\ku\R} \simeq \underline{\prod_{p\ge 0}\mathbf{H}\R[2p]}\simeq  \prod_{p\ge 0} \mathbf{H}\Omega_{\R}[2p]\to \DD^{\per}\ ,$$ 
where the last map is the natural inclusion. This composition provides the lower horizontal map in the pull-back square  in $\Sh_{\Sp}(\Mf_{\C\textrm{-}\fol})$ \begin{equation}\label{wefwefewfewfewfewfw}
\xymatrix{\mathbf{H}_{\R,\mathrm{Del}}\ar[d]\ar[r]&\DD^{-}\ar[d]\\
\underline{\ku\R}\ar[r]&\DD^{\per}}
\end{equation}
which defines $$\mathbf{H}_{\R,\mathrm{Del}}\in \Sh^{h}_{\Sp}(\Mf_{\C\textrm{-}\fol})\ .$$ On the one hand, this  is the Hodge-filtered version of $\ku\R$-theory. On the other hand,  it is a generalization  of  real Deligne cohomology to foliated manifolds. In fact, for a smooth complex algebraic variety $X$ we have a natural isomorphism \begin{equation}\label{frff877823r87786r8723r32r32r}
\pi_{*}(\mathbf{H}_{\R,\mathrm{Del}} (X^{\an},T^{0,1}X))\cong \prod_{p\in \nat} H_{\mathrm{Del},\mathrm{an}}^{2p-*}(X^{\an},\R(p))\ .
\end{equation}
Note that $\mathbf{H}_{\R,\mathrm{Del}} $   does not  involve the weight-filtration and therefore   reflects  the  ``wrong" Hodge filtration on $H^{*}(X^{\an};\C)$ for non-proper $X$.
 
 \bigskip
 
The natural map $\ku\to \ku\R$ induces a morphism of pull-back squares $\eqref{ffewfewfwefewfewf234}\to \eqref{wefwefewfewfewfewfw}$ and therefore a morphism,
$$\ch_{\R,\mathrm{Del}}:  \ku^{\flat}\to\mathbf{H}_{\R,\mathrm{Del}}\ .$$
The composition
$$\bK\stackrel{\reg}{\to} \ku^{\flat}\stackrel{\ch_{\R,\mathrm{Del}}}{\to} \mathbf{H}_{\R,\mathrm{Del}}$$ yields indeed  Beilinson's regulator if one applies this to the  foliated manifolds $(X^{\an},T^{0,1}X)$, precomposes with
$\bK^{\alg}(X)\to \bK(X^{\an},T^{0,1}X)$, see Example \ref{kdjqwkldjwqldjwqldwqd}, and uses the identification \eqref{frff877823r87786r8723r32r32r}. 
This easily follows from the description of Beilinson's regulator given in  \cite{Bunke:2012fk}, \cite{Bunke:2013aa}.
}
\end{rem}




Let $(M,\cF)$ be a foliated manifold.

\begin{lem}
If $\codim(\cF)<p$, then
we have a natural isomorphism
$$\ku\C/\Z^{-p-1}(M)\cong  ku^{\flat,-p}(M,\cF)\ .$$
\end{lem} \proof This easily follows from
 $$\pi_{p}(\DD^{-}(M,\cF))\cong 0\cong   \pi_{p+1}(\DD^{-}(M,\cF))\ .$$ \hB

\begin{kor}
If $\codim(\cF)<p$, then the regulator (Definition \ref{ojlfewfefwefw234234}) induces a map
$$\widetilde \reg:K^{-p}(M,\cF)\to \ku\C/\Z^{-p-1}(M)\ .$$
\end{kor}

%
%

 \begin{rem}\label{2ddhi3dhio32doi2doo2oidjud2}{\rm 
 We have a factorization of $\reg_{X}$ defined in Definition \ref{klfwefewfewfwf}
 as
 $$\xymatrix{K_{X}^{-p}(*)\ar[r]^{\reg_{X}}\ar[d]^{\eqref{dkejwdlewdewded}}&\ku\C/\Z^{-p-1}(X)\\
 K^{-p}(X,\cF_{\max})\ar[ur]^{\widetilde \reg}\ar[r]^{\eqref{cejkwdjejdhewdhewkdei832e32e32e32e2e}}_{\cong}&\bK(\C)^{-p}(X)\ar[u]^{\sigma}}$$
 for every $p\ge 1$. Here $\sigma:\bK(\C)\to \ku\C/\Z$ is the morphism discussed e.g. in \cite[7.21]{karoubiast}, \cite[Ex. 6.9]{2013arXiv1311.3188B}.}
\end{rem}


\begin{rem}{\rm  
In \cite{Bunke:2014aa} we asked whether the map  \begin{equation}\label{hdgqhjdgqwhdghjqwgdwhjgdj7861} 
K_{p}(C^{\infty}(X))\to K^{\mathrm{top}}_{p}(C^{\infty}(X))\end{equation} can be non-trivial for $p>\dim(X)$.
 This question has an analogue in the foliated case.

Note that
$$\kappa^{*}\bK\to \underline{\ku}$$ (see Example \ref{kjffjewfewjkfhewkjfewfkewfhewf87z} for $\kappa$) is the homotopification morphism. 
  The question is now: \begin{prob}\label{lkdeqwdqwdqwd} Let $(M,\cF)$ be a foliated manifold and $p\in \nat$ be such that $\codim(\cF)<p$.   Is the map  
$$K^{-p}(M,\cF)\to \ku^{-p}(M)$$
 trivial?
 \end{prob}

In the special case of a minimal foliation we ask wether $$K^{-p}(X,\cF_{\min})\to \ku^{-p}(X)$$
can be non-trivial for $p>\dim(X)$.
 The difference to \eqref{hdgqhjdgqwhdghjqwgdwhjgdj7861} can best be explained by the commutative diagram  
$$ \xymatrix{K_{p}(C^{\infty}(X))\ar[d]\ar[r]&K_{p}^{\mathrm{top}}(C^{\infty}(X))\ar[d] \\K^{-p}(X,\cF_{\min})\ar[r]&\ku^{-p}(X)}\ ,$$
where the   vertical maps  are induced by   sheafification in the $X$-direction.
  
 \bigskip
 
In the foliation case   we can answer the Question \ref{lkdeqwdqwdqwd} affirmatively at least rationally. 
\begin{prop} 
Let $(M,\cF)$ be a foliated manifold, $p\in \nat$ such that $\codim(\cF)<p$ and
$x\in K^{-p}(M,\cF)$.
 Then the image
$x_{\Q}\in \ku\Q^{-p}(M)$ vanishes.\end{prop}\proof
Note that the natural map $\ku\Q^{*}(M)\to \ku\C^{*}(M)$ is injective and that
the Bockstein sequence
$$\ku\C/\Z^{-p-1}(M)\stackrel{\beta}{\to} \ku^{-p}(M)\stackrel{c}{\to} \ku\C^{-p}(M)$$
is exact.  We write $x_{\C}$ for the image of $x_{\Q}$ in $\ku\C^{-p}(M)$.
We have
$x_{\C}=c(\beta(\tilde\reg(x))=0$. \hB

}\end{rem}

\bibliographystyle{alpha}
\bibliography{rho}
\end{document}